\documentclass[review,hidelinks,onefignum,onetabnum]{siamart220329}

\usepackage{mathrsfs}%
\usepackage{algorithm}%
\usepackage{verbatim}

\newtheorem{example}[theorem]{Example}
\newtheorem{assumption}{Assumption}
\newtheorem{Remark}[theorem]{Remark}

\newcommand{\Rb}{\mathbb{R}}
\newcommand{\Sb}{\mathbb{S}}

\newcommand{\Jcal}{\mathcal{J}}
\newcommand{\Ical}{\mathcal{I}}
\newcommand{\Mcal}{\mathcal{M}}
\newcommand{\Ncal}{\mathcal{N}}
\newcommand{\Scal}{\mathcal{S}}
\newcommand{\Lcal}{\mathcal{L}}
\newcommand{\Vcal}{\mathcal{V}}
\newcommand{\sign}{\text{sign}}

\usepackage{subcaption}
\usepackage{graphicx}

\usepackage{lipsum}
\usepackage{amsfonts}
\usepackage{graphicx}
\usepackage{epstopdf}
\usepackage{algorithmic}
\ifpdf
  \DeclareGraphicsExtensions{.eps,.pdf,.png,.jpg}
\else
  \DeclareGraphicsExtensions{.eps}
\fi

\newsiamremark{remark}{Remark}
\newsiamremark{hypothesis}{Hypothesis}
\crefname{hypothesis}{Hypothesis}{Hypotheses}
\newsiamthm{claim}{Claim}

\headers{AVOIDING STRICT SADDLE POINTS OF NONCONVEX
REGULARIZED PROBLEMS}{Luwei Bai, Yaohua Hu, Hao Wang, Xiaoqi Yang}

\title{AVOIDING STRICT SADDLE POINTS OF NONCONVEX REGULARIZED PROBLEMS\thanks{Submitted to the editors DATE.
\funding{This work was funded by the Fog Research Institute....}}}

\author{Luwei Bai\thanks{School of Information Science and Technology, ShanghaiTech University, Shanghai, P. R. China 
  (\email{bailuwei.kexiny@gmail.com}).}
\and Yaohua Hu\thanks{School of Mathematical Sciences, Shenzhen University, Shenzhen, P. R. China 
  (\email{mayh-hu@szu.edu.cn}).}
\and  Hao Wang\thanks{Corresponding author. School of Information Science and Technology, ShanghaiTech University, Shanghai, P. R. China 
  (\email{wanghao1@shanghaitech.edu.cn}).}
\and Xiaoqi Yang\thanks{Department of Applied Mathematics, The Hong Kong Polytechnic University, Hong Kong, P. R. China 
  (\email{mayangxq@polyu.edu.hk}).}}

\usepackage{amsopn}

\ifpdf
\hypersetup{
  pdftitle={Avoiding strict saddle points of nonconvex regularized problems},
  pdfauthor={Luwei Bai, Yaohua Hu, Hao Wang, Xiaoqi Yang}
}
\fi

\begin{document}\nolinenumbers

\title{Avoiding strict saddle points of nonconvex regularized problems}
\maketitle

\begin{abstract}
This paper considers a class of nonconvex and nonsmooth sparse optimization problems, encompassing most existing nonconvex sparsity-inducing terms.  We show that their second-order optimality conditions depend only on the nonzeros of the stationary points.  We propose two damped iteratively reweighted algorithms, which are the iteratively reweighted $\ell_1$ algorithm (DIRL$_1$) and the iteratively reweighted $\ell_2$  (DIRL$_2$)  algorithm,  to solve these problems.   For DIRL$_1$, we show that the reweighted $\ell_1$ subproblem has the support identification property so that DIRL$_1$ locally reverts to a gradient descent algorithm around a stationary point. 
For DIRL$_2$, we show that the solution map of the reweighted $\ell_2$ subproblem is differentiable and Lipschitz continuous everywhere. Therefore, the solution maps of DIRL$_1$ and DIRL$_2$ and their inverses are Lipschitz continuous, and the strict saddle points are their unstable fixed points. By applying the stable manifold theorem, these algorithms starting from almost every initial point are shown to converge to local minima when the strict saddle point property is assumed.
\end{abstract}

\begin{keywords}
strict saddle property, iterative reweighted $\ell_1$ algorithm, nonconvex regularization, center stable manifold theorem
\end{keywords}

\begin{AMS}  49J50, 49M20, 49M37, 65K05, 65K10, 90C26, 90C30
\end{AMS}

\section{Introduction}\label{sec1}

In this paper, we consider the following nonconvex and nonsmooth regularization problem,
\begin{equation}
    \min\limits_{x\in\mathbb{R}^n} F(x) :=  f(x) + \lambda\psi(x), \tag{P} \label{pp} 
\end{equation}
where $f:\mathbb{R}^n \rightarrow \mathbb{R}$ is a twice continuously differentiable function and  $\psi(x)$ is a penalty/regularization term taking the form $\psi(x) := \sum\limits_{i=1}^n s(x_i)$ with $s(x_i) = r(|x_i|)$ and $r$ satisfying the following properties.
\begin{assumption}\label{ass.1}
 $r$ is continuous and concave on $[0,+\infty)$ with $r(0) = 0$ and C$^2$ smooth on $(0,+\infty)$, $\lim\limits_{t\rightarrow 0^+}r'(t) > 0$ and $r'(t)\ge 0$ for all $t > 0$.
\end{assumption}
\noindent This assumption indicates $r(t) > 0$ for all $t >0$, but $r$ could be nonconvex or not Lipschitz continuous around 0. 

The regularization term $\psi$ in \eqref{pp} encourages sparse solutions, offering a favorable balance between model complexity, interpretability, and predictive performance of the model. This makes it well suited for a wide range of tasks including machine learning, statistics, signal processing, and image analysis \cite{figueiredo2007majorization, nikolova2008efficient, saab2008stable, yin2015minimization}. Many widely used nonconvex regularizations, such as the EXP penalty  \cite{bradley1998feature}, the LOG penalty \cite{lobo2007portfolio}, and the $\ell_p$ quasi-norm penalty ($0<p<1$) \cite{frank1993statistical}, satisfy Assumption \ref{ass.1} and
can offer advantages including improved sparsity and better model interpretability. 
 We provide the explicit form of several commonly used nonconvex regularizations in Table \ref{table1}.   It is worth noting that these regularizers can be categorized into two classes: (i) Lipschitz continuous regularizers with a bounded subdifferential around 0 ($r'(|0^+|) < \infty$); (ii)  nonLipschitz continuous regularizers with unbounded subdifferential around 0 ($r'(|0^+|) = \infty$) such as the $\ell_p$ quasi-norm (LPN).

Numerous studies focus on gradient-based optimization methods that are guaranteed to find first-order stationary points. 
However, the nonconvexity of $F$ can lead to multiple local minima and saddle points, which can cause algorithms to converge to a stationary point rather than a local minimum.  In addition, the presence of saddle points can significantly slow the convergence of optimization algorithms, especially gradient-based methods. In this situation the gradient provides misleading information about the descent directions at a saddle point.  As a result, avoiding saddle points can help optimization algorithms more effectively converge to better solutions.  Therefore, it becomes a prevalent topic, particularly in high-dimensional optimization problems, and is also the central focus of this paper.

 In smooth optimization, many methods have been developed and analyzed to escape saddle points, including first-order methods \cite{panageas2019first, li2019alternating, sun2019escaping}, momentum-based methods \cite{wang2021escaping}, second-order algorithms \cite{paternain2019newton}, the algorithmic frameworks for constrained nonlinear problems \cite{mokhtari2018escaping}, stochastic algorithms \cite{ge2015escaping, daneshmand2018escaping}, and techniques explicitly designed to identify and circumvent saddle points during optimization, such as the random perturbation technique \cite{jin2017escape, jin2018accelerated, lu2019pa, criscitiello2019efficiently}.  
In smooth settings, it is widely assumed that \emph{strict saddle property} holds, which means that all saddle points have at least one strictly negative eigenvalue. This property intuitively implies that at any saddle point, there exists a direction in which the objective function strictly decreases, and it has been shown to be generally satisfied for smooth problems \cite{bhojanapalli2016global,sun2018geometric}.

In nonsmooth optimization, avoiding saddle points poses a greater challenge. 
At stationary points, traditional tools such as the (sub)differential may not reliably indicate descent directions or facilitate avoidance of saddle points. Developing nonsmooth optimization algorithms being capable of effectively converging to desirable solutions and circumventing saddle points is inherently more complex than that in smooth optimization. 
 An example of this complexity was demonstrated in the analysis of the Douglas splitting algorithm for nonsmooth problems, as discussed in \cite{atenas2024weakly}. Recent research  \cite{davis2022proximal}  has made strides in defining saddle points in nonsmooth cases as those occurring on the (smooth) active manifold by assuming that $0\in \text{rint} \partial F(x)$ is satisfied at the stationary point. Furthermore, several proximal algorithms with damped methods have been proposed to converge to local minima in nonsmooth optimization settings, including the proximal point algorithm, the proximal gradient algorithm, and the damped proximal linear algorithm \cite{davis2022proximal}. Furthermore, an inexact analog of a stochastically perturbed gradient method applied to the Moreau envelope has also been introduced \cite{davis2022escaping}. 
 In these studies, the nonsmooth function $r(x)$ is assumed to be $\rho$ weakly convex, which means that $r(\cdot) + \frac{\rho}{2}\|\cdot\|^2$
is convex.  
They assumed an active manifold around the limit point and then defined the strict saddle points as those of the $f$ restricted to the manifold.   The weak $\rho$-convexity of $r$ plays a crucial role in constructing a $C^2$ differentiable Moreau envelope with the same stationary points and strict saddle points in the proximal methods.

As for the regularized problem \eqref{pp}, several first-order methods have been proposed. \cite{xu2012l_}  introduced a proximal gradient method tailored for solving the $\ell_p$ regularized problem with $0<p<1$. This method was subsequently extended to address the $\ell_{p,q}$ group norm regularized problems with $0<q<1\leq p$ in \cite{hu2017group}. Furthermore, the difference of convex (DC) algorithm was enlisted \cite{ahn2017difference} to handle nonconvex regularization problems, given that many existing nonconvex regularizers $r(x)$ can be reformulated as DC functions.
It is important to note, however, that the DC algorithm generally assumes $r'(|0^+|)<\infty$ to accommodate the DC decomposition, which does not hold for the $\ell_p$ regularized problem. Another avenue of research involves iteratively reweighted algorithms, which approximate 
the nonconvex regularizers by iteratively reweighted $\ell_1$ or $\ell_2$ regularization subproblems. Initially applied to solve the unconstrained $\ell_p$ regularized problem in \cite{Chartrand2008Iteratively} without convergence analysis, this approach was later employed to address general constrained nonconvex regularization problems in \cite{wang2021nonconvex}. Subsequent work in \cite{lu2014iterative} established its global convergence and demonstrates its locally stable support property \cite{wang2021relating}. Furthermore, convergence rate analysis under the Kurdyka-\L{}ojasiewicz property \cite{wang2023convergence} and its extrapolated version was proposed in \cite{wang2022extrapolated}. 
While these first-order methods primarily focus on identifying first-order stationary points, none of them guarantee the avoidance of strict saddle points.

\subsection{The main results}

In this paper, our main purpose is to develop and analyze first-order methods that can provably avoid strict saddle points of the nonconvex regularization problem \eqref{pp} without the assumption of weak convexity.  In particular, we focus on iteratively reweighted methods which approximate the regularization term by simpler weighted $\ell_1$ or $\ell_2$ norms \cite{chen2013optimality, chen2010convergence, lai2011unconstrained, lu2014iterative,sun2017global, wang2021nonconvex,wang2021relating,wang2022extrapolated}. In contrast to the proximal gradient algorithm, the iteratively reweighted $\ell_1$ method (${\rm IRL}_1$)   solves a sequence of weighted $\ell_1$ subproblems, which have an explicit solution characterized by the soft-threshold operator. Therefore, it  is widely applicable to various nonconvex regularization problems. 
This method offers a favorable local property \cite{wang2021relating}, known as the model/support identification, under the following additional assumption. 
\begin{assumption} \label{ass.2}
 For any stationary point $x^*$,  $0 \in {\rm rint}\partial F(x^*)$ (where $\partial F$ denotes the regular subdifferential and ${\rm rint}$ denotes the relative interior.)  \end{assumption}  
\noindent This kind of assumption is also needed to guarantee the existence of an active manifold around a stationary point as assumed by \cite{davis2022proximal}. 
  At the tail end of the ${\rm IRL}_1$ iteration, the sign of the iterates remains unchanged while the nonzero components are kept uniformly bounded away from 0.   In other words, the algorithm identifies a model subspace defined by an index set
$\mathcal{I}$ as follows after a finite number of iterations. 
\begin{definition} 
The model subspace defined by  $\mathcal{I} \subset \{1,\ldots, n\}$  is  
$$\mathcal{M}_\mathcal{I}:= \{u\in\mathbb{R}^n :  u_j = 0, \forall j\in \mathcal{I}^c\}.$$
\end{definition}
\noindent This means that every limit point of the iterates is also contained in $\mathcal{M}_\mathcal{I}$. 
   This property motivates us to treat the mapping of  ${\rm IRL}_1$  as smooth (locally around a limit point) in $\mathcal{M}_\mathcal{I}$, although the underlying optimization problem is nonsmooth in $\Rb^n$.

The iteratively reweighted $\ell_2$ method (${\rm IRL}_2$) is also widely applied to these nonconvex regularization problems,  which solves a sequence of weighted $\ell_2$ subproblems and the analytic solution of each subproblem is trivially given. Compared to ${\rm IRL}_1$, ${\rm IRL}_2$ does not enjoy the model subspace identification/active-manifold identity property. In fact, it can be often witnessed that this algorithm can generate iterates with elements asymptotically approaching zero.  
The behavior of ${\rm IRL}_2$ brings another challenge in analyzing the avoidance of strict saddle points. Current analysis for nonsmooth optimization algorithms ubiquitously depends on the active-manifold identification property, which ensures that the algorithm reverts to a smooth algorithm on the active manifold in a local neighborhood of the stationary point.  However,  there is no answer to how a nonsmooth first-order algorithm like ${\rm IRL}_2$ would behave in the local region of a strict saddle point.

We now summarize our contributions as follows: 

\begin{enumerate} 
\item[(i)]  We show that the local minimum and second-order optimality of $F$ are equivalent to those of $F$ restricted in $\mathcal{M}_\mathcal{I}$, and we define the strict saddle points of $F$ as those in $\mathcal{M}_\mathcal{I}$.  

\item[(ii)]  We develop damped IRL$_1$ algorithms individually for the two types of regularizers, $r'(|0^+|)<\infty$ and $r'(|0^+|)=\infty$, meaning that the next iterate is the combination of a traditional IRL$_1$ iterate and its previous iterate. 
The only difference between these two cases is that a relaxation $\epsilon>0$ is added to the
regularizers for the latter case to prevent $r'(\cdot)$ from becoming $\infty$ at each iteration.   At the end of the iterations, we show that this algorithm reverts to a damped gradient descent algorithm in the subspace $\mathcal{M}_\mathcal{I}$ and reduces the rest components to 0 at a linear speed.   We then apply the center stable manifold theorem and prove that the algorithm can avoid convergence to strict saddle points.

\item[(iii)]  We develop damped IRL$_2$ algorithms for nonconvex and nonsmooth regularized problems and provide global convergence.  
Under an additional assumption which is satisfied by the $\ell_p$ regularization, we show that the algorithm is $C^1$ smooth and nondegenerate 
everywhere.  We also prove the equivalence between a unstable fixed point of the iteration map and a strict saddle point of the nonconvex regularized problem.  Therefore, the avoidance of
strict saddle points is derived.  This is the first algorithm shown to have this property without using the active-manifold identification property or the stochastic technique in nonsmooth optimization setting. 

\end{enumerate}


\subsection{Outline}
The remainder of this paper is organized as follows. \S 2 introduces preliminary results and the definition of the active strict saddle point of problem $\eqref{pp}$. 
\S 3 proposes the damped iteratively reweighted $\ell_1$ algorithms as well as its global convergence and the property of avoiding strict saddle points.  \S 4 provides the damped iteratively reweighted $\ell_2$ algorithms, its global convergence, and the property of avoiding the strict saddle points.

\subsection{Notation}
We consider $\mathbb{R}^n$ as the real $n$-dimensional Euclidean space with inner product $\langle\cdot, \cdot\rangle$ and norm $\|x\| = \sqrt{\langle x, x\rangle}$. For $x\in\mathbb{R}^n$, $x_i$ denotes the $i$-th element of $x$. Let $\bar{\mathbb{R}}=(-\infty, \infty]$. Furthermore, $\mathbb{R}_{+}^n$ denotes the nonnegative orthant of $\mathbb{R}^n$ with $\mathbb{R}^n_{+} = \{x\in\mathbb{R}^n : x\geq 0\}$, $\mathbb{R}^n_{-} = \{x\in\mathbb{R}^n : x\leq0\}$, and so forth for $\bar{\mathbb{R}}^n_+ = \{x\in\bar{\mathbb{R}}^n : x_i \in[0,+\infty], i=1,...,n\}$. Let  $I$ and  $\boldsymbol{0}$ be the identity matrix and zero matrix of appropriate sizes, respectively. 
Given a symmetric matrix $H\in\mathbb{R}^{n\times n}$, let $\lambda_{\min}(H)$ be 
the smallest eigenvalue of $H$. 

We denote $\|\cdot\|_p$ as the quasi-norm $\ell_p$ with $p\in(0, 1)$, that is, $\|x\|_p = (\sum_{i=1}^n |x_i|^p)^{\frac{1}{p}}$. This does not define a proper norm because of the lack of subadditivity. The index sets of nonzero and zero components of $x$ are defined as $\mathcal{I}(x) = \{i \mid x_i \neq 0\}$, and $\mathcal{J}(x) = \{i \mid x_i = 0\}$, respectively. Moreover, $\mathcal{I}^*$ and $\Jcal^*$ denotes $\mathcal{I}(x^*)$ and $\mathcal{J}(x^*)$ for short. 
${\rm diag}(x)$ is the diagonal matrix by placing  elements of vector $x$ on the main diagonal. 
${\rm sign}(x) : \mathbb{R}^n \rightarrow \mathbb{R}^n$ is the sign function defined by $[{\rm sign}(x)]_i=1$ if $x_i > 0$, $[{\rm sign}(x)]_i=-1$ if 
$x_i < 0$ and $[{\rm sign}(x)]_i = 0$ if $x_i =0$. 

$DT(x)$ denotes the Jacobian matrix of a map $T : \mathbb{R}^n \rightarrow \mathbb{R}^n$ at $x$.
Given $x\in \Rb^n$ and $\Ical\subset \{1,\ldots, n\}$,  denote  $x_{\mathcal{I}}  \subset \mathbb{R}^{\mathcal{I}}  $ or $[x_i]_{i\in\Ical}$  as the subvector of the elements in  
$ \Ical$; also,  
given  $f : \mathbb{R}^n \rightarrow \mathbb{R}$ smooth with respect to $x_i, i\in \Ical  $, 
denote $\nabla_{\mathcal{I}} f(x)$ and $\nabla^2_{\mathcal{II}} f(x)$  as  the 
gradient and the Hessian of $f$ at $x$ with respect to $\Ical$, respectively. 
If $f$ is smooth, $\nabla^2_{\mathcal{II}} f(x) =[ \nabla^2 f(x)]_{\mathcal{II}}$ and   also let 
$\nabla^2_{\mathcal{IJ}} f(x) =  [ \nabla^2 f(x)]_{\mathcal{IJ}}$ for $J\subset \{1,\ldots, n\}$. Denote $f_{| \Ical}$ as the restriction of $f$ in   $\Mcal_{\Ical}$, i.e., by fixing 
$x_i\equiv 0, i\in \{1,\ldots, n\}\setminus \Ical$; therefore, 
$f_{| \Ical}(x_{\Ical})  := f([x_{\mathcal{I}}; 0_{\Jcal}]).$  
 Let $x\in\mathbb{R}^n$ and 
$w\in \mathbb{R}_{+}^n$, the soft-thresholding operator $\mathbb{S}_w: \mathbb{R}^n\to \mathbb{R}^n $ with   $w$  defined as 
$ [\mathbb{S}_{ w }(x)]_i = \text{sign}(x_i)  \max(|x_i|-w_i,0)$. 
In particular, if $w_i = \infty$,  it follows that $[\mathbb{S}_{ w }(x)]_i = 0$. 
%

\begin{table}[t]  
    \renewcommand{\arraystretch}{2}
    \centering
    \begin{tabular}{c|c|c|c|c}
        \hline
        Regularizer & $\psi(x)$ & $r(|x_i|)$ & $r'(|x_i|)$ & $r''(|x_i|)$ \\
        \hline
        EXP \cite{bradley1998feature}& $\sum\limits_{i=1}^n(1 - e^{-p|x_i|})$ & $1 - e^{-p|x_i|}$ & $pe^{-p|x_i|}$ & $-p^2e^{-p|x_i|}$\\
        \hline
        LOG \cite{lobo2007portfolio} & $\sum\limits_{i=1}^n \log(1+p|x_i|)$ & $\log(1+p|x_i|)$ &  $\frac{p}{1+p|x_i|}$ & $-\frac{p^2}{(1 +p|x_i| )^2}$ \\
        \hline
        FRA \cite{fazel2003log} & $\sum\limits_{i=1}^n \frac{|x_i|}{|x_i| + p}$ & $\frac{|x_i|}{|x_i| + p}$ & $\frac{p}{(|x_i| + p)^2}$ & $-\frac{2p}{(|x_i| + p)^3}$  \\
        \hline
        LPN \cite{frank1993statistical} & $\sum\limits_{i=1}^n |x_i|^p$ & $ |x_i|^p$ & $p|x_i|^{p-1}$ & $p(p-1)|x_i|^{p-2}$  \\ 
        \hline
        TAN \cite{candes2008enhancing} & $\sum\limits_{i=1}^n \arctan(\frac{|x_i|}{p})$ & $\arctan(\frac{|x_i|}{p})$ &$\frac{p}{|x_i|^2 + p^2}$ & $-\frac{2p|x_i|}{(|x_i|^2 + p^2)^2}$  \\
        \hline
    \end{tabular}
    \centering
    \caption{ Concrete examples of nonconvex regularization }
    \label{table1}
\end{table}

\section{Optimality and strict saddle points}

\subsection{Optimality conditions} We discuss the first- and second-order optimality conditions of problem \eqref{pp} in this subsection. First, let us recall the definition of subgradients and the first-order optimality condition.
\begin{definition}(subgradients \cite{rockafellar2009variational})
    Consider a function $f:\mathbb{R}^n \rightarrow \mathbb{\bar{R}}$ and $f(\bar{x})$ is finite. For any $v\in\mathbb{R}^n$, one says that 
    \begin{enumerate}
        \item[(1)]  $v$ is a regular subgradient of $f$ at $\bar{x}$, if $v\in\hat{\partial}f(\bar{x})$ where
            \begin{equation}
                \hat\partial f(\bar{x}) = \{v \mid f(\bar{x}) + \langle v, x-\bar{x}\rangle + o(\|x - \bar{x}\|) \leq f(x),\  {\rm as}\ x\rightarrow \bar{x}\}.
            \end{equation}
        \item[(2)]  $v$ is a (general) subgradient of $f$ at $\bar{x}$, written $v\in\partial f(\bar{x})$, if there are sequences $x^k \xrightarrow[f]{} \bar{x}$ and $v^k\in\hat{\partial} f(x^k)$ with $v^k\rightarrow v$.
    \end{enumerate}
\end{definition}

 It is trivial to see
$s'(0^+) = r'(|0^+|)$ and $s'(0^-) = - r'(|0^+|)$. 
The subgradients of $s$ are given in the following proposition. 
\begin{proposition}\label{subgradient}  
It holds that $\partial F(x) = \hat\partial F(x) = \nabla f(x)+ \lambda \partial \psi(x)$ where 
 \begin{equation}\label{psi.sub}
  \partial \psi(x) =   \hat\partial \psi(x) = \partial \sum\limits_{i=1}^n s( x_i ) 
                   =   \partial s( x_1 ) \times ... \times \partial s( x_n )
  \end{equation}
and  
 \begin{equation}\label{r.sub}
       \partial s( x_i ) = \hat\partial s( x_i )  =  \begin{cases}
               \{ \text{sign}(x_i) r'(|x_i|)\},         &i\in\Ical(x),\\
              [- r'(|0^+|) , r'(|0^+|)],   &i\in\Jcal(x) \text{ and }  r'(|0^+|) < \infty, \\
            \mathbb{R},                  &i\in\Jcal(x) \text{ and } r'(|0^+|) = \infty. 
        \end{cases}      
         \end{equation}
\end{proposition}

\begin{proof}
We only need to calculate the subgradients of $s$ at $0$, i.e., $i\in \Jcal(x)$.  
If $r'(|0^+|)<\infty$, define for any $x\in\mathbb{R}$, 
\[ \hat s(x) :=  \begin{cases}   r(|x|)  & \text{   if  }  x \ge 0,\\
r'(|0^+|)x & \text{ otherwise}.
\end{cases}\]
Obviously, $\hat s(\cdot)$ is smooth on $\mathbb{R}$ with $\hat s'(0) = r'(|0^+|) > 0$ and
$s(x) = \max\{ \hat s(x), - \hat s(x) \}$. 
By \cite[Exercise 8.31]{rockafellar2009variational}, we have $\partial s(0) = [- r'(|0^+|) ,   r'(|0^+|)]$. 
On the other hand,  $\liminf\limits_{x \rightarrow 0, x \neq 0} \frac{s(x) - s(0) - v (x - 0)}{|x - 0|} \geq 0$ holds   for  $v = s'(0^+) = r'(|0^+|)$ and $ v = s'(0^-) =  - r'(|0^+|)$, which means $r'(|0^+|)\in \hat\partial s(0)$ and $ -r'(|0^+|)\in \hat\partial s(0)$. 
  It follows from \cite[Theorem 8.6]{rockafellar2009variational} that 
  $  \partial s(0)  = [ -r'(|0^+|) , r'(|0^+|)]    \subset \hat\partial s(0)  \subset   \partial s(0)$. 
  Hence,    $ \partial s(0) =\hat\partial s(0)   = [ -r'(|0^+|) , r'(|0^+|)]$.

If $r'(|0^+|) = \infty$, we have $\liminf\limits_{x\rightarrow 0, x\neq 0} \frac{r(|x^*|) - r(0) - v(x - 0)}{|x - 0|} \geq 0$ for all $v\in\mathbb{R}$. Therefore, we have $\mathbb{R}\subseteq\hat{\partial} s(0) \subseteq \partial s(0) \subseteq \mathbb{R}$. Hence, $\hat{\partial} s(0) = \partial s(0) = \mathbb{R}$.

In general, we have shown that \eqref{r.sub} is true.   It then follows from \cite[Proposition 10.5]{rockafellar2009variational}  that  \eqref{psi.sub} is true. 
\end{proof}

We can appeal to the classical \emph{rule of Fermat} \cite[Theorem 10.1]{rockafellar2009variational} to characterize a local minimum. 
If $x^*$ is a local minimum of \eqref{pp}, then
    \begin{align}
        0 \in \partial F(x^*) = \nabla f(x^*) + \lambda\partial\psi(x^*)  \label{opt.cond0}
    \end{align}
  This condition,  combined with Assumption \ref{ass.2} and the explicit form of  the subgradients, 
yields the following results. 

\begin{proposition}\label{Firstoptimal2}
   If $x^*$ is a local minimum  for \eqref{pp}, then  the following hold
\begin{align} 
            &\nabla_i f(x^*) = -\lambda\text{sign}(x_i^*) r'(|x^*_i|),   &&i\in\mathcal{I}^*, \label{opt.cond1}\\ 
            &\nabla_i f(x^*) \in (-\lambda  r'(|0^+|), +\lambda  r'(|0^+|)),   &&i\in\Jcal^*.  \label{opt.cond2}
 \end{align}
\end{proposition} 
\noindent Notice that if $r'(|0^+|) = \infty$, \eqref{opt.cond2} holds naturally true.   Any point that satisfies the necessary optimality conditions \eqref{opt.cond1} and \eqref{opt.cond2} is referred to as a stationary point.  


 In the following, we demonstrate that the second-order optimality condition can be characterized by nonzero components.  
  

\begin{theorem} \label{the2ndCon} 
\begin{enumerate}
\item[(i)] (necessity) 
     If   $x^*$  is a    local minimum  of $F$, then   
       \begin{align}
        d^T\nabla^2f(x^*)d + \lambda\sum\limits_{i\in \mathcal{I}^*}r''(|x^*_i|) d_i^2 \ge & \  0, \quad \forall d\in \mathcal{M}_{\Ical^*}, \label{eq2ndCon}
    \end{align} 
    or, equivalently,  $\lambda_{\min}(\nabla^2_{\mathcal{I}^*\Ical^*} F(x^*)) \ge 0$. 
\item[(ii)] (sufficiency)   A stationary point $x^*$ of $F$ is a strict  local minimum   if 
    \begin{align}
        d^T\nabla^2f(x^*)d + \lambda\sum\limits_{i\in \mathcal{I}^*}r''(|x^*_i|) d_i^2 > & \  0, \quad \forall d\in \mathcal{M}_{\Ical^*}\setminus\{0\},  \label{eq2ndCon2}
    \end{align} 
    or, equivalently,  $\lambda_{\min}(\nabla^2_{\mathcal{I}^*\Ical^*} F(x^*)) > 0$. 
    \end{enumerate}

\end{theorem} 
\begin{proof} (i) is proved in \cite[Theorem 3.3]{chen2013optimality}.  We only prove (ii). 

    Assume by contradiction that $x^*$ is not a strict local minimum. Then there exists a sequence of nonzero vectors $\{u^k\}$ satisfying \begin{equation}
       F(x^* + u^k) \le  F(x^*)\quad {\rm and}\quad \lim\limits_{k\rightarrow\infty}u^k = 0. \label{ThSecond1}
    \end{equation}
      Consider $F_{|\Ical^*}$, 
    which is locally smooth around $x^*_{\mathcal{I}^*}$ and 
    $x^*_{\mathcal{I}^*}$ is stationary for $F_{|\Ical^*}$. 
    It follows from $u^k = \begin{bmatrix}
        u_{\Ical^*}^k ; u_{\Jcal^*}^k  
    \end{bmatrix}$ that 
    \begin{equation}
        \label{check.2}
        \begin{aligned}
            &F(x^*+u^k) - F(x^*)\\ 
                                   = &\ f(x^*+u^k) + \lambda\sum\limits_{i\in\mathcal{I}^*}r(|x^*_i + u_i^k|) + \lambda\sum\limits_{i\in\Jcal^*}r(|x^*_i + u_i^k|)-F_{|\Ical^*}(x^*_{\mathcal{I}^*})  \\
                                   = & \ f(x^*+u^k) + \lambda\sum\limits_{i\in\mathcal{I}^*}r(|x^*_i + u_i^k|) + f(x^*+u_{\mathcal{I}^*}^k) \\
                                   & - f(x^*+[u_{\mathcal{I}^*}^k;\boldsymbol{0}]) + \lambda\sum\limits_{i\in\Jcal^*}r(|u_i^k|)-F_{|\Ical^*}(x^*_{\mathcal{I}^*}) \\
                                   = &\ f(x^*+u^k) - f(x^*+[u_{\mathcal{I}^*}^k;\boldsymbol{0}]) + F_{|\Ical^*}(x^*_{\mathcal{I}^*}+u_{\mathcal{I}^*}^k) - F_{|\Ical^*}(x^*_{\mathcal{I}^*}) +  \lambda\sum\limits_{i\in\Jcal^*}^k r(|u_i^k|).
        \end{aligned}
        \end{equation}
    From the classic optimality condition for a smooth function and \eqref{eq2ndCon2}, it is easy to see that $x^*$ is a strict local minimum of $F_{|\Ical^*}(x)$. Therefore for sufficiently large $k$, 
    \begin{equation}
        F_{|\Ical^*}(x^*_{\mathcal{I}^*}+u_{\mathcal{I}^*}^k) - F_{|\Ical^*}(x^*_{\mathcal{I}^*}) > 0. \label{ThSecond4}
    \end{equation}
    Moreover, by the Lagrangian mean value theorem,  
    \begin{equation}\label{check.1} r(|u_i^k|)  =  r(|u_i^k|) - r(0) =     r'(\hat u_i^k) | u^k_i|, \text{  \  with  }  \hat u_i^k \in (0, |u_i^k|). 
    \end{equation}
    The concavity (Assumption \ref{ass.1}) of $r$ implies that $ r'(\hat u_i^k) \ge r'(\| u^k_{\Jcal^*}\|_\infty)$, which, together  with \eqref{check.1}, yields that   
\[    \begin{aligned}
     \sum\limits_{i\in\Jcal^*}r(|u_i^k|)   \ge    \sum\limits_{i\in\Jcal^*} 
   r'(\| \hat u^k_{\Jcal^*}\|_\infty)  \|u^k_{\Jcal^*}\|_1  \ge \sum\limits_{i\in\Jcal^*} 
   r'(\| \hat u^k_{\Jcal^*}\|_\infty)  \|u^k_{\Jcal^*}\|_2.  
    \end{aligned} 
  \]
    This, combined with \eqref{ThSecond4} and \eqref{check.2}, yields
    \begin{align}\label{final.ineq}
        F(x^* + u^k) - F(x^*) > (\lambda r'(\|\hat u^k_{\Jcal^*}\|_\infty) + \eta^k)\|u^k_{\Jcal^*}\|,
    \end{align}
    where 
    \[  \eta^k:= \frac{   f(x^*+u^k) - f(x^*+[u_{\mathcal{I}^*}^k;\boldsymbol{0}])  }{ \|u_{\Jcal^*}^k\| } =  \nabla f(x^*+[u_{\mathcal{I}^*}^k;\boldsymbol{0}]+ t [\boldsymbol{0};u_{\mathcal{J}^*}^k])^T   [\boldsymbol{0};u_{\mathcal{J}^*}^k)]/\|u_{\Jcal^*}^k\|,\] 
    and $t \in [0,1]$, 
    implying $$\limsup\limits_{k\to\infty} |\eta^k| \le \limsup\limits_{k\to\infty}\|\nabla f(x^*+[u_{\mathcal{I}^*}^k;\boldsymbol{0}]+ t [\boldsymbol{0};u_{\mathcal{J}^*}^k])\|_\infty\| \tfrac{u_{\mathcal{J}^*}^k}{\|u_{\Jcal^*}^k\| }\|_1 = \|\nabla f(x^*)\|_\infty.$$ 
    It follows from  
     $  \| \nabla f(x^*)\|_\infty < \lambda r'(|0^+|)$ by  \eqref{opt.cond2} that      
    $|\eta^k| < \lambda r'(\|\hat u^k\|_\infty) $ for sufficiently large $k$,  as  
    $ \|\hat u^k\|_\infty \to 0$. 
    It follows from  \eqref{final.ineq}
   that for sufficiently large $k$, 
    \[         F(x^* + u^k) - F(x^*) \geq (\lambda r'(\|\hat u^k\|_\infty) + \eta^k)\|u^k_{\Jcal^*}\| > 0,\]
contradicting  \eqref{ThSecond1}.
\end{proof}

\subsection{  Strict Saddle Points}
In this subsection, we define the strict saddle point of the problem \eqref{pp}. In smooth settings, the strict saddle point can be defined by virtue of a negative minimum eigenvalue of the corresponding Hessian, indicating a descent direction. 
For our case,  Theorem \ref{the2ndCon} implies that a saddle point can be defined using nonzeros.  
\begin{definition}
    A stationary point $x^*$ of $F(x)$ is a  strict saddle point if 
    \[ d^T \nabla^2_{\Ical^*\Ical^*} F(x^*) d < 0\] for $\forall d\in\mathbb{R}^{|\Ical^*|}\setminus\{0\}$. If all stationary points of $F(x)$ are strict saddle points or local minima, then we say that $F$ satisfies the strict saddle property. \label{defSaddle}
\end{definition}

This definition is related to the ``active strict saddle point'' proposed in  \cite{davis2022proximal},  in which it is shown that the presence of negative curvature alone in nonsmooth settings does not guarantee escape from saddle points. To address this, they define the \emph{active strict saddle point} for nonsmooth problems as the saddle point on the active manifold containing the
stationary point on which the function is partly smooth.  Our definition also belongs to this kind, since $\mathcal{M}_{\Ical^*}$ is
an active manifold that contains $x^*$.

To show that the strict saddle property is mathematically generic,
many works show that this property holds for a full measure of linear perturbation \cite[Theorem 2.9]{davis2022proximal}, \cite[Theorem 2]{bianchi2023stochastic}.  Namely,  the perturbed function 
 has the strict saddle property for almost all linear perturbations. 
 In our case, we also have the same results by considering the following linearly perturbed nonsmooth regularization problem for $v\in\mathbb{R}^n$,
\begin{equation}
\min\limits_{x\in\mathbb{R}^n} F_v(x) := F(x) - \langle v, x\rangle =  f(x) + \lambda\psi(x) - \langle v, x\rangle.  \label{ppv} 
\end{equation}
  To show the desired result,  we need the assistance of the famous Sard's theorem  \cite{matsumoto2002introduction}. 
 Our proof of the strict saddle property mainly follows the proof of \cite[Lemma 2.7]{ubl2023linear}, which is first proved in \cite{milnor2025lectures}.


\begin{definition} Let $f : \Mcal \rightarrow \mathcal{N}$ be a smooth map of manifolds. If the derivative at $x$ is singular (Jacobian is \emph{not} of full rank), then $x$ is called a critical point of $f$.
\end{definition} 

\begin{theorem}(Sard's Theorem)
Let $f: \Mcal \rightarrow \mathcal{N}$ be a smooth map of manifolds, and $C$ be the set of critical points of $f$ in $\Mcal$. Then $f(C)$ has measure zero in $\Ncal$.
\end{theorem}

The genericity of the strict saddle property is stated in the following. 
\begin{theorem} Let $\Ical$ stand for $\Ical(x)$. 
    For a full Lebesgue measure set of perturbations $v\in\mathbb{R}^n$, the perturbation function $F_v(x) := F(x) -  \langle v, x\rangle$ satisfies
    \begin{equation}
    \inf\limits_{d\in \mathbb{R}^{|\Ical|} \setminus\{0\}}d^T \nabla^2_{\Ical\Ical}F(x)d \neq 0, \label{st.saddle}
    \end{equation}
    for all stationary points $x$ of $F_v$. \label{lem.generic} 
\end{theorem} 

\begin{proof} Let $s \in \{+1, - 1, 0\}^n$ be a vector of sign pattern and 
\[ \Mcal (s) = \{  x  \in \Rb^n   \mid \sign(x) = s\}  \]
 be a manifold.  Consider the behavior of $F_v$ in $\Mcal(s)$, where 
 $\Ical = \Ical(x) = \Ical(s)$ for all $x\in \Mcal(s)$.  Let $v:  \Mcal(s)  \rightarrow \Rb^{n}$ be a smooth map: 
 $ [v(x)]_{\Ical}  = \nabla_{\Ical} F(x) $ and $[v(x)]_{\Ical^c} = x$. 
The Jacobian of $v$ at $x$ is then given by 
$ Dv(x) = \begin{bmatrix}
  \nabla_{\Ical\Ical}^2 F(x) & 0 \\ 0 & I  
\end{bmatrix}$, 
and $Dv(x)$ has full rank if and only if $\nabla_{\Ical\Ical}^2 F(x)$ has full rank. 
By Sard's Theorem, $\{ v(x)  \mid \nabla_{\Ical\Ical}^2 F(x) \text{  is singular}, x\in \Mcal(s) \}  \text{ has measure 0 in } \mathbb{R}^n$,
which implies its subset 
\[ \Vcal(s):= \{ v \mid 0\in \partial F_v(x),  \inf\limits_{d\in \mathbb{R}^{|\Ical(s)|} \setminus\{0\}} d^T \nabla^2_{\Ical\Ical}F(x)d = 0, x \in \Mcal(s) \}\]
 has   measure  0 in $\mathbb{R}^n$.
    
 Now we consider all the possible sign pattern vectors $s \in \{+1, - 1, 0\}^n$. There are at most $3^n$ different patterns.  Therefore,  
\[ \bigcup_{s\in \{+1, - 1, 0\}^n}   \Vcal(s)  = \{ v \mid 0\in \partial F_v(x), \inf\limits_{d\in \mathbb{R}^{|\Ical|} \setminus\{0\}}d^T \nabla^2_{\Ical\Ical}F(x)d = 0, x \in \mathbb{R}^n \}\]   has   measure  0 in $\mathbb{R}^n$.
\end{proof}

\subsection{Center Stable Manifold Theorem} \label{subsecCSM}

As with many other works in smooth settings like \cite{lee2019first}, we will interpret our algorithm as the fixed point iteration of a nonlinear map $T$. 
Likewise, we define fixed points with at least one eigenvalue of $DT$ larger than $1$ in magnitude as \emph{unstable fixed points}, as described below.

\begin{definition} (Unstable fixed point)
    Let $x^*$ be a fixed point on a $C^1$ map $T: \mathbb{R}^{n} \rightarrow \mathbb{R}^{n}$.   
    $x^* $ is called an unstable fixed point of $T$ if  $DT(x^* )$ has at least one eigenvalue of magnitude strictly larger than 1.  \label{def.fix.point}
\end{definition}

The Center Stable Manifold Theorem (CSM) plays a crucial role in the avoidance of strict saddle points. 
The CSM theorem establishes a favorable property for the fixed-point iteration of a $C^1$ local diffeomorphism map: the set of initial points around a fixed point, from which an algorithm converges to the fixed point with at least one eigenvalue of $DT$ greater than $1$ in magnitude, has Lebesgue measure zero.

\begin{theorem} (The Center Stable Manifold Theorem)
    Let $0$ be a fixed point for the $C^1$ local diffeomorphism $T: U \rightarrow \mathbb{R}^n$ where $U$ is a neighborhood of zero in $\mathbb{R}^n$. Let $E^s \oplus E^c \oplus E^u$ be the invariant split of $\mathbb{R}^n$ into the generalized eigenspaces of $DT(0)$ corresponding to eigenvalues of absolute value less than one, equal to one and greater than one. There exists a $C^1$ embedded disk $W_{loc}^{cs}$ that is tangent to $E_s \oplus E_c$ at 0 and a neighborhood $B$ around 0 such that $T(W_{loc}^{cs}) \cap B \subset W_{loc}^{cs}$.  In addition,  if $T^k(x) \in B$ for all $k\ge 0$, 
    then $x  \in   W_{loc}^{cs}$. \label{the2}
\end{theorem}

If $D T(0)$ is invertible,  then $T$ is locally a diffeomorphism around 0. 
The CSM theorem only guarantees the avoidance of strict saddle points in the neighborhood for local diffeomorphism,  and this local result needs to be extended to the global setting to guarantee that the set of almost all initial points from which an algorithm  converges to any unstable fixed point has a measure zero.  This requires that the map has the so-called 0-property, meaning the preimage of a measure zero set under this map still has measure zero. 
This result was established in \cite[Theorem 2]{lee2019first} for a global diffeomorphism map and then for a global lipeomorphism (A map $T: \mathbb{R}^n \rightarrow \mathbb{R}^n$ is called a lipeomorphism if it is a bijection and both $T$ and $T^{-1}$ are Lipschitz continuous) in \cite{davis2022proximal}. 
 Moreover,  to extend this property from a particular stationary point $x^*$ to the set of all stationary points,  this can be shown using the fact that $\Rb^n$ is second countable no matter if the number of stationary points is countable (see \cite[Corollary 2.12]{davis2022proximal} for a standard proof).  
  We list their results below. 
  
\begin{theorem}\cite[Theorem 2]{lee2019first}
    Let $g$ be a $C^1$ mapping from $\mathbb{R}^n \rightarrow \mathbb{R}^n$ and $det(DT(x)) \neq 0$ for all $x\in\mathbb{R}^n$.  Let $\mathcal{U}_T$ consist of all unstable fixed points $x$ of $T$
    Then the set of initial points from which an algorithm  converges to an unstable fixed point has measure zero, 
    $\mathscr{U} = \{ x\in \mathbb{R}^n : \lim_{k\to \infty} T^k(x) \in \mathcal{U}_T\}$ has zero Lebesgue measure. 
\end{theorem}
 
%

\begin{corollary} \cite[Corollary 12]{davis2022proximal} \label{coro.saddle}
Let $T: \mathbb{R}^n\to \mathbb{R}^n$ be a lipeomorphism and let $\mathcal{U}_T$ consists of 
all unstable fixed points $x$ of $T$ at which the Jacobian $D  T(x)$ is invertible. Then, the set of initial 
points attracted by such fixed points 
$ \mathscr{U} = \{ x\in \mathbb{R}^n : \lim\limits_{k\to \infty} T^k(x) \in \mathcal{U}_T\}$
has zero Lebesgue measure. 
\end{corollary}

To summarize, the analysis of  escaping from strict saddle points under the strict saddle property generally consists of three steps: 
\begin{enumerate} 
\item[(1)] Interpret the iterative algorithm as a fixed point iteration and build the equivalence of the unstable fixed points and the strict saddle points. 
\item[(2)] Show that the fixed point iteration is locally a $C^1$  local diffeomorphism around a strict saddle point (e.g., a $C^1$ smooth map with invertible Jacobian at the fixed point) so that the CSM theorem guarantees the avoidance of converging to a strict saddle point. 
\item[(3a)] Show that the fixed-point iteration is globally a diffeomorphism.
\item[(3b)]  (Alternative to (3a)) Show that the fixed-point iteration is globally a lipeomorphism. 
\end{enumerate}

\section{Damped Iterative Reweighted $\ell_1$ Algorithms} \label{sec3}
The ${\rm IRL}_1$ algorithm is an effective method to solve problem \eqref{pp}.   
At each iteration, the ${\rm IRL}_1$ algorithm solves the following weighted $\ell_1$ regularized subproblem,  
\begin{align*}
 x^{k+1} =   S^x(x^k)  = \arg\min_{y} G(y;x^k):=\nabla f(x^k)^T(y - x^k) + \tfrac{\beta}{2}\|y - x^k\|^2 + \lambda\sum_{i=1}^n \tilde w^k_i |y_i|,
\end{align*}
where weight $\tilde w_i^k = r'(|x_i^k|)$.  

 
This algorithm may quickly converge to local minima due to the rapid growth of the weight values, especially for the case where $r'(|0^+|) =  \infty$ and $w(0) = \infty$.  To address this issue, a common approach is to add a   perturbation 
$\epsilon_i  > 0$ to each $r$ and drive $\epsilon_i$ to 0 during the iteration, resulting in the following approximation $F(x, \epsilon)$ of $F(x)$ 
\begin{equation}
    \min\limits_{x\in\mathbb{R}^{n}, \epsilon\in\mathbb{R}^{n}_+} F(x,\epsilon) :=  f(x) + \lambda\sum\limits_{i=1}^{n} r(|x_i|+\epsilon_i). \label{probepsi} 
\end{equation}
Here we   treat $\epsilon$ also as variables,  so that the algorithm updates $x$ and $\epsilon$ by 
\[ (x^{k+1}; \epsilon^{k+1}) = S(x^k,\epsilon^k),\] 
where $S=[S^x; S^\epsilon]$ and 
\begin{align*}
 x^{k+1} &=    S^x(x^k,\epsilon^k) = \arg\min\limits_{y} G(y;x^k,\epsilon^k): = \nabla f(x^k)^Ty + \tfrac{\beta}{2}\|x^k - y\|^2 + \lambda\sum\limits_{i=1}^n w_i^k |y_i|, \\
 \epsilon^{k+1} &=  S^\epsilon(x^k, \epsilon^k) =   \mu\epsilon^k
\end{align*}
with $w_i^k := w(x_i^k, \epsilon_i^k) = r'(|x_i^k| + \epsilon_i^k)$   
and     $\mu\in(0,1)$.

The subproblem solution can be given using the soft-thresholding operator, 
\begin{equation}\label{operator.soft}
\begin{aligned}
    S^x(x^k,\epsilon) = \mathbb{S}_{\lambda w^k/\beta} (x^k - \tfrac{1}{\beta} \nabla f(x^k)), 
\end{aligned}
\end{equation}
or explicitly, 
\begin{equation}\label{operator.s}
   S^x_i(x^k,\epsilon^k) = \begin{cases} x_i^k - \tfrac{1}{\beta}(\nabla_i f(x^k ) + \lambda w^k_i ), &  \text{if }  x_i^k - \tfrac{1}{\beta}(\nabla_i f(x^k) + \lambda w^k_i ) > 0,
    \\   x_i^k - \tfrac{1}{\beta}(\nabla_i f(x^k) - \lambda w^k_i ), &  \text{if } x_i^k - \tfrac{1}{\beta}(\nabla_i f(x^k) - \lambda w^k_i ) < 0,
    \\   0,  &  \text{otherwise }. \end{cases}  
\end{equation}

To take advantage of the CSM theorem, the map $S$ needs to be a local diffeomorphism at a stationary point and a lipeomorphism globally. However, this is not true in our case, since $S$ is generally noninvertible. 
In particular, $S$ is not a bijection, as it maps multiple values to $0$.
  Motivated by the technique used in \cite{davis2022proximal},  
  we propose to use the  damped operator: $T=  (1-\alpha) I  + \alpha S$, 
%
which can be written explicitly as 
\begin{equation}\label{fix.point} 
(x^{k+1}, \epsilon^{k+1}) = T(x^k,\epsilon^k) = \begin{bmatrix}
    (1-\alpha) x^k + \alpha S^x(x^k, \epsilon^k) \\
    (1-\alpha(1-\mu))\epsilon^k
\end{bmatrix}
\end{equation}
with $\epsilon_0 \in \Rb^n_{++}$ and $\alpha \in (0, 1)$. 

 This algorithm, hereinafter nicknamed DIRL$_1$, is presented in Algorithm \ref{alg1}.  If the algorithm does not involve the vector of $\epsilon$, we can simply treat this case as setting $\epsilon^0 = 0$ in  Algorithm \ref{alg1} and the iteration simply reverts to 
  \[ x^{k+1} = T(x^k)   = (1-\alpha)x^k + \alpha S(x^k).\]

We have the following result to connect the optimality of \eqref{pp} and the fixed point of our iteration. 
\begin{proposition}\label{prop.fix.point} 
 $x^*$ is a stationary point of \eqref{pp} if and only if $x^* = S^x(x^*,0)$. In addition,  
 $x^*$ is a stationary point of \eqref{pp} if and only if $(x^*,\epsilon^*) = T(x^*,\epsilon^*)$. 
\end{proposition} 
\begin{proof} Suppose $x^*$ is a stationary point. It follows from \eqref{opt.cond1} that  if $x_i^* > 0$, 
$x_i^* -  \frac{1}{\beta}[\nabla_if(x^*) + \lambda w(x_i^*, 0)]  = x_i^* - 
\frac{1}{\beta}[\nabla_if(x^*) + \lambda \text{sign}(x_i^*) r'(|x_i^*|)] =  x_i^* > 0$, implying 
$S_i^x(x^*,0)  = x_i^*$.  This result also holds  by the same argument if $x_i^* < 0$.  

As for $x_i^* = 0$,   $x_i^* -  \frac{1}{\beta}[\nabla_if(x^*) + \lambda w(x_i^*, 0)]  = - \frac{1}{\beta}[\nabla_if(x^*) + \lambda r'(|0^+|)] < 0 $ 
by  \eqref{opt.cond2}.  The same argument also yields 
$x_i^* -  \frac{1}{\beta}[\nabla_if(x^*) - \lambda w(x_i^*, 0)]  = - \frac{1}{\beta}[\nabla_if(x^*) - \lambda r'(|0^+|)]   > 0$. 
Therefore, $S_i^x(x^*,0) = 0 = x_i^*$.  
In general, we have shown that $x^*=  S^x(x^*,0)$ in all cases.

Conversely,  suppose $x^*= S^x(x^*, 0)$. 
If $x_i^* > 0$, it must be true that $x_i^* = x_i^* -  \frac{1}{\beta}[\nabla_if(x^*) + \lambda w(x_i^*, 0)]    = x_i^* - \frac{1}{\beta}[\nabla_if(x^*) + \lambda r'(|x_i^*|)] > 0 $, implying 
$ 0 = x_i^* - S_i^x(x^*,0) = - \frac{1}{\beta}[\nabla_if(x^*) + \lambda r'(|x_i^*|)] $. 
The same argument yields $ 0 = x_i^* - S_i^x(x^*,0) = - \frac{1}{\beta}[\nabla_if(x^*) - \lambda r'(|x_i^*|)]$ if 
$x_i^*<0$. Therefore,  
$\nabla_if(x^*) + \lambda \text{sign}(x_i^*)  r'(|x_i^*|) = 0$ for $i\in \Ical^*$. 

 If $x_i^* = 0$, it indicates that 
 $| \frac{1}{\beta} \nabla_if(x^*) |  =  | x_i^* -  \frac{1}{\beta} \nabla_if(x^*) | \le \tfrac{\lambda}{\beta} w(0, 0)    = \tfrac{\lambda}{\beta} r'(|0^+|)$. 
Hence, $\nabla_i f(x^*) \in [- \lambda r'(|0^+|), - \lambda r'(|0^+|)]$. 
In general, $x^*$ satisfies \eqref{opt.cond0} and is a stationary point of \eqref{pp}. 
By Assumption \ref{ass.2},  it satisfies \eqref{opt.cond1} and \eqref{opt.cond2}. 

The rest part of the statement of this proposition is obviously true by noticing 
that $(x^*;\epsilon^*) = T(x^*,\epsilon^*)$ implies that $\epsilon^* = 0$. 
\end{proof}
%

\begin{algorithm}
    \begin{algorithmic}[1] 
        \caption{Damped Iteratively Reweighted $\ell_1$ Algorithm (${\rm DIRL}_1$)}
        \label{alg1}
        \STATE \textbf{Input} $\beta > \alpha L_{\nabla f}/2$, $\alpha\in(0, 1)$, $x^0$, $\epsilon^0 > 0$ and $\mu \in (0, 1)$.
        \STATE \textbf{Initialize}: set $k = 0$.
        \REPEAT 
            \STATE \textrm{Reweighting:} $w(x_i^k, \epsilon_i^k) = r'(|x_i^k| + \epsilon_i^k)$.
            \STATE \textrm{Compute new iterate:} 
            \STATE $y^{k} = \arg\min\limits_{y\in\mathbb{R}^n} \{\nabla f(x^k)^T(y - x^k) + \frac{\beta}{2}\|y - x^k\|^2 + \lambda \sum\limits_{i=1}^{n}w(x_i^k, \epsilon_i^k)|y_i|\}$.
            \STATE $x^{k+1} = (1-\alpha)x^k + \alpha y^{k}$.
            \STATE $\epsilon^{k+1} = (1-\alpha(1-\mu))\epsilon^k$.
            \STATE \textrm{Set} $k = k + 1$.
        \UNTIL{Convergence}
    \end{algorithmic}  
\end{algorithm}

\subsection{Convergence analysis}

It should be noted that the same conclusion for Algorithm \ref{alg1} can be readily deduced by setting $\epsilon$ to zero. 
We first make the following assumptions for the algorithm.
\begin{assumption}
Given the initial point $(x^0$, $\epsilon^0)$ and $F^0 = F(x^0, \epsilon^0)$, the level set 
$\mathcal{L}(F^0) := \{x \mid F(x ) \leq F^0\}$
is bounded. \label{ass.3}
\end{assumption}

\noindent 
This assumption means that there exists $C > 0$ such that 
 \begin{equation}\label{xnablabd}  \| x- \tfrac{1}{\beta} \nabla f(x)\|_{\infty} \leq C,  \quad \forall x \in \Lcal(F^0). \end{equation}
It also means $f$ is Lipschitz continuous and differentiable with a Lipschitz gradient on $\mathcal{L}(F^0)$. Hence, there exist constants $L_{f}$ and $L_{\nabla f}$, such that $\forall  x,y\in \Lcal(F^0)$,
 \begin{equation}\label{f.nonexp}  
  |f(x)-f(y)| \le L_f\|x-y\| \text{ and }       \|\nabla f(x) - \nabla f(y)\| \leq L_{\nabla f}\|x - y\|. 
     \end{equation}
Moreover, \eqref{f.nonexp} implies that  if $r'(|0^+|) = \infty$,  then
for any stationary point $x^*\in \Lcal(F^0)$, 
$|x_i^*| \ge (r')^{-1}(L_f/\lambda), \forall i\in\Ical^*$ by \eqref{opt.cond1} and 
$r''(|x_i^*|) > 0$ is bounded above. As for $r'(|0^+|) < \infty$, it holds naturally that 
$r''(|x_i^*|) $ is bounded above for $i\in\Ical^*$. 

Overall,  we can conclude that for any stationary point $x^* \in \mathcal{L}(F^0)$, there exists $\rho>0$ such that  
\begin{equation}\label{hess is bounded}
\| \nabla^2 _{\Ical^*\Ical^*} F(x^*) \|
 \le \rho. 
\end{equation}

We first show the monotonicity of $F(x, \epsilon)$.  Considering  the model $G$ at $(x^k,\epsilon^k)$ 
\begin{equation*}
    G(x; x^k, \epsilon^k) := \nabla f(x^k)^Tx + \frac{\beta}{2}\|x^k - x\|^2 + \lambda\sum_{i=1}^n w(x_i^k, \epsilon_i^k) |x_i|\
\end{equation*}
with minimum $y^k$, its reduction 
caused by  $x^{k+1} = (1-\alpha)x^k + \alpha y^{k}$  is defined as 
\[     \Delta G(x^{k+1}; x^k, \epsilon^k) = \ G(x^k; x^k, \epsilon^k) - G(x^{k+1}; x^k, \epsilon^k).\]
The reduction of   $F$ caused by $(x^{k+1},\epsilon^{k+1})$ with $\epsilon^{k+1}  = (1-\alpha(1-\mu))\epsilon^k$ 
  is defined as 
\begin{equation*}
\begin{aligned}
    \Delta F(x^{k+1}, \epsilon^{k+1}) = &\ F(x^k, \epsilon^k) - F(x^{k+1}, \epsilon^{k+1}). 
\end{aligned}
\end{equation*}

\begin{proposition}
    Suppose Assumptions \ref{ass.1}-\ref{ass.3} hold. Let $\{(x^k, \epsilon^k)\}$ be the sequence generated by Algorithm \ref{alg1} with $\alpha \in (0,1)$ and  $\beta \geq \frac{\alpha L_{\nabla f}}{2}$. It holds that \{$F(x^k, \epsilon^k)$\} is monotonically decreasing and that the reduction satisfies
    \begin{equation*}
        F(x^0, \epsilon^0) - F(x^k, \epsilon^k) \geq (\frac{\beta}{\alpha} - \frac{L_{\nabla f}}{2})\sum\limits_{t=0}^{k-1}\|x^t - x^{t+1}\|^2, 
    \end{equation*}
    which means  
     $ \lim\limits_{k\rightarrow\infty} \|x^{k+1} - x^{k}\| = 0$. Moreover,    $\{ x^k \}\subset \mathcal{L}(F^0)$ and \eqref{xnablabd} holds 
    for any $x = x^k, k \in \mathbb{N}$. \label{prop}
\end{proposition}
\begin{proof}
     Since the concavity of $r(\cdot)$ on $\mathbb{R}_{++}$ gives $r(x_1) \leq r(x_2) + r'(x_2)(x_1 - x_2)$ for any $x_1, x_2 \in \mathbb{R}_{++}$, we have that for $i = 1, ..., n$,
    \begin{equation*}
        r(|x_i^{k+1}|+\epsilon_i^k) \leq r(|x_i^k|+\epsilon_i^k) + r'(|x_i^k|+\epsilon^k_i)(|x_i^{k+1}| - |x_i^k|).
    \end{equation*}
    Summing up the above inequality over $i$ yields
    \begin{equation}
        \sum\limits_{i=1}^{n}r(|x_i^{k+1}|+\epsilon_i^k) - \sum\limits_{i=1}^{n}r(|x_i^k|+\epsilon_i^k)  \leq  \sum\limits_{i=1}^{n}w(x_i^k, \epsilon^k_i)(|x_i^{k+1}| - |x_i^{k}|). \label{p2.con2}
    \end{equation}
On the other hand, it follows that 
    \begin{equation}
        f(x^k) - f(x^{k+1}) \geq \nabla f(x^k)^T(x^k - x^{k+1}) - \frac{L_{\nabla f}}{2}\|x^k - x^{k+1}\|^2. \label{p2.con1}
    \end{equation} 
Combining \eqref{p2.con2} and \eqref{p2.con1}, we have 
\begin{equation}
\begin{aligned}
 &  \Delta F(x^{k+1}, \epsilon^{k+1}) \geq 
 F(x^k, \epsilon^{k}) - F(x^{k+1}, \epsilon^k) \\
 \geq \ &  \nabla f(x^k)^T(x^k - x^{k+1})  - \frac{L_{\nabla f}}{2}\|x^k - x^{k+1}\|^2 + \lambda\sum\limits_{i=1}^{n}w(x_i^k, \epsilon_i^k)(|x_i^k| - |x_i^{k+1}|) \\ 
  = \ & \Delta G(x^{k+1}; x^k, \epsilon^k) + \frac{\beta - L_{\nabla f}}{2}\|x^k - x^{k+1}\|^2,
  \end{aligned}  \label{p2.con4}
  \end{equation}
    where the first inequality is due to the fact that $\epsilon^{k+1} \leq \epsilon^k$. 
    
    For each subproblem, $y^k$ satisfies the optimality condition
    \begin{equation}
        \nabla f(x^k) + \lambda W(x^k, \epsilon^k)\xi^k + \beta(y^k - x^k)= \boldsymbol{0}.  \label{p2.opt}
    \end{equation}
    with $\xi^k \in \partial\|y^k\|_1$ and $W(x^k, \epsilon^k) := {\rm diag}(w(x_1^k, \epsilon_1^k), ..., w(x_n^k, \epsilon_n^k))$.
    Consequently, by $x^{k+1} = (1-\alpha)x^k + \alpha y^k$ and the optimality condition \eqref{p2.opt}, the following holds 
    \begin{equation} \label{p1.con5}
    \begin{aligned}
     & \ \Delta G(x^{k+1}; x^k, \epsilon^k) \\
    =&\ \nabla f(x^k)^T(x^k - x^{k+1}) + \lambda\sum\limits_{i=1}^{n}w(x_i^k,\epsilon_i^k)(|x_i^k| - |x_i^{k+1}|) - \frac{\beta}{2}\|x^k - x^{k+1}\|^2 \\
        \ge &\ \nabla f(x^k)^T(\alpha x^k - \alpha y^k) + \lambda\sum\limits_{i=1}^{n}w(x_i^k,\epsilon_i^k)(\alpha|x_i^k| - \alpha|y_i^k|) - \frac{\beta\alpha^2}{2}\|x^k - y^k\|^2 \\
        \geq&\ \alpha [ \nabla f(x^k)^T(x^k - y^k) + \lambda\sum\limits_{i=1}^{n}w(x_i^k, \epsilon_i^k)\xi_i^{k}(x_i^k - y_i^k)] - \frac{\beta\alpha^2}{2}\|x^k - y^k\|^2 \\
        =&\ \alpha[\nabla f(y^k) + \lambda W(x^k, \epsilon^k)\xi^{k} + \beta(y^k - x^k)]^T(x^k - y^k)  + \frac{\alpha\beta(2 - \alpha)}{2}\|x^k - y^k\|^2 \\
        =&\ (\frac{\beta}{\alpha} - \frac{\beta}{2})\|x^k - x^{k+1}\|^2, 
    \end{aligned}
    \end{equation} 
            where the first inequality is by $\alpha \in (0,1)$ and the second is by the convexity of $|\cdot|$.  
    Combining \eqref{p2.con4} and \eqref{p1.con5}, we get
    \begin{equation*}
        \Delta F(x^{k+1}, \epsilon^{k+1}) \geq (\frac{\beta}{\alpha} - \frac{L_{\nabla f}}{2} )\|x^k - x^{k+1}\|^2.
    \end{equation*}
    Replacing $k$ with $t$ and summing up from $t = 0$ to $k-1$, we have
    \begin{equation*}
        \sum\limits_{t=0}^{k-1}[F(x^t, \epsilon^t) - F(x^{t+1}, \epsilon^{t+1})] = F(x^0, \epsilon^0) - F(x^k, \epsilon^k) \geq (\frac{\beta}{\alpha} - \frac{L_{\nabla f}}{2})\sum\limits_{t=0}^{k-1}\|x^t - x^{t+1}\|^2.
    \end{equation*}
    It follows that $\{(x^k, \epsilon^k)\} \in \mathcal{L}(F^0)$ since $F(x^k, \epsilon^k) \leq F(x^0, \epsilon^0)$. By Assumption \ref{ass.3}, it is bounded 
    and $\mathop{\lim\inf}\limits_{k \rightarrow \infty} F(x^k, \epsilon^k) > -\infty$. Taking $k \rightarrow \infty$, we have  
    \begin{equation*}
        (\frac{\beta}{\alpha} - \frac{L_{\nabla f}}{2})\sum\limits_{t=0}^{\infty}\|x^t - x^{t+1}\|^2 \leq F(x^0, \epsilon^0) - \mathop{\lim\inf}\limits_{k \rightarrow \infty} F(x^k, \epsilon^k) < \infty.
    \end{equation*}
    which implies $\lim\limits_{k \rightarrow \infty}\|x^k - x^{k+1}\|^2 = 0$.

    Moreover, $F(x^k) \le F(x^k,\epsilon^k) \le F(x^0,\epsilon^0)$, meaning $\{x^k\}\subset \mathcal{L}(F^0)$. Therefore,  \eqref{xnablabd} holds for all $x = x^k, k \in \mathbb{N}$.    
\end{proof}

With the monotonicity of $F$, we derive the global convergence of the ${\rm DIRL}_1$ method. 
\begin{theorem}(Global convergence)\label{GlobalConverge.p1}
    Suppose Assumptions \ref{ass.1} and \ref{ass.3} hold. Let the sequence $\{(x^k,\epsilon^k)\}$ be generated by ${\rm DIRL}_1$ and $\Omega$ be the set of the limit points of $\{x^k\}$. Then $\epsilon^*= 0$ and $\Omega \neq \emptyset$ and any $x^*\in \Omega$ is a stationary point of \eqref{pp}.  
\end{theorem}
\begin{proof}
Firstly, the boundedness of $\mathcal{L}(F^0)$ from Assumption \ref{ass.3} implies $\Omega \neq \emptyset$. Then let $x^*$ be a limit point of $\{x^k\}$ with subsequence $\{x^k\}_{\mathcal{S}} \rightarrow x^*$ and $\{y^k\}_{\Scal} \rightarrow x^*$ since $y^k = \frac{1}{\alpha}(x^k - x^{k+1}) + x^k$. This also implies   $\sign(y_i^k) = \sign(x_i^*), i\in \Ical^*$ for sufficiently large $k \in \Scal$.   
The optimality condition of the $k$th subproblem \eqref{p2.opt} implies 
$ - \beta(y_i^k - x_i^k) \in \nabla_if(x^{k})  + \lambda w(x_i^k, \epsilon_i^k)\partial |y_i^k|.
$
    If $w(x_i^*, 0) < +\infty$,  it follows from the  outer semicontinuity of 
      $\partial r$
that  there exists $\xi^*_i\in\partial |x^*_i|$ such that
$
      0  = \nabla_i f(x^*) + \lambda w(x_i^*, 0)  \xi^*_i%
$
    as $y_i^k - x_i^k \overset{\Scal}{\to} 0$, $  x_i^k \overset{\Scal}{\to} x_i^*$ and 
 $\epsilon_i^k \to 0$. 
Therefore, for $w(x_i^*,0) = r'(|x_i^*|) < \infty, i\in \Ical^*$, \eqref{opt.cond1} is satisfied;    
for  $w(x_i^*,0) = r'(|0^+|) = \infty, i\in \Jcal^*$,   \eqref{opt.cond2} is satisfied by Assumption \ref{ass.2}.   Therefore, $x^*$ is a stationary point.
  \end{proof}

\subsection{Avoidance of strict saddle points}\label{sec3.1}

The map $S^x$ has the property of so-called ``active manifold identification'', 
which means that the iterates fall onto a smooth manifold after finite iterations.  The iteratively reweighted $\ell_1$ algorithms were proved to have this property by \cite{wang2022extrapolated, wang2023convergence} and the manifold in this case is simply the subspace consisting of nonzeros at a limit point. 


\begin{theorem} (Model Identification and Local Smoothness) \label{ModelId}    
    Under Assumptions \ref{ass.1}-\ref{ass.3}, for any stationary point $(x^*,\epsilon^*)$, there exists a sufficiently small neighborhood $B$ of $(x^*,\epsilon^*)$, so that the following hold: 
    \begin{enumerate}
    \item[(i)]  For any point $(x,\epsilon)\in B$,   $S^x(x, \epsilon) \in (x^*+\mathcal{M}_{\Ical^*})$ and $ \text{sign}(S^x(x, \epsilon)) = \text{sign}(x) = \text{sign}(x^*)$.
    \item[(ii)] $S$ is $C^1$ smooth in $B$ and $S^x(x,\epsilon) = x - \frac{1}{\beta} g(x,\epsilon)$ with 
    \begin{equation}
         g_i(x,\epsilon) = \left\{\begin{aligned}
    &\beta x_i,  & i \in \Jcal^*,\\
    &\nabla_i f(x) + \lambda r'(x_i\text{sign}(x_i^*) + \epsilon_i)\text{sign}(x_i^*), & i\in\Ical^*.
\end{aligned}\right. \label{smoothprob.grad}
    \end{equation}
    \item[(iii)] Map $T$ is $C^1$-smooth in $B$.  
        \end{enumerate}
\end{theorem}
\begin{proof}
 (i) Since $S$ is a continuous map and $x^* = S^x(x^*,\epsilon^*)$, we have \[ \text{sign}(x_i)= \text{sign}(x^*_i) \text{\ and\ }  \text{sign}(S_i(x, \epsilon)) =  \text{sign}(S_i(x^*,0))   = \text{sign}(x^*_i),  \quad i\in\Ical^*\] for any $(x,\epsilon)$ sufficiently close to $(x^*,\epsilon^*)$.   It remains to show that for any $i\in\Jcal^*$, $S_i(x, \epsilon) = 0$ for any $(x,\epsilon)\in B$. 
By Assumption \ref{ass.2}, $  | x_i^* - \tfrac{1}{\beta} \nabla_i f(x^*) | = |\tfrac{1}{\beta} \nabla_i f(x^*)| < \tfrac{1}{\beta} \lambda w(0,0) $. 
    Therefore,  there exists a neighborhood $B$ such that  $  | x_i - \tfrac{1}{\beta} \nabla_i f(x) |  <\tfrac{1}{\beta} \lambda w(x_i, \epsilon_i) $, $\forall (x,\epsilon) \in B$, which is equivalent to  $    x_i - \tfrac{1}{\beta} ( \nabla_i f(x) - \lambda w(x_i, \epsilon_i))  < 0  $ 
    and $   x_i - \tfrac{1}{\beta} (  \nabla_i f(x) -   \lambda w(x_i, \epsilon_i) )   > 0$, meaning $S_i(x,\epsilon) = 0$. 

(ii)     It follows from (i) that there exists a neighborhood $B$ around $(x^*,\epsilon^*)$, so that the subproblem at $(x^k,\epsilon^k)\in B$ agrees with the following constrained problem,
\begin{equation}\label{sub.prob.eqv}
    \begin{aligned}
        S^x(x^k,\epsilon^k) = & \  \arg\min\limits_{x\in x^*+\mathcal{M}_{\Ical^*}}\{\nabla f(x^k)^Tx + \tfrac{\beta}{2}\|x-x^k\|^2 + \lambda\sum_{i=1}^n w^k_i \text{sign}(x_i^*) x_i \}\\
         = & \  \arg\min\limits_{x\in \mathbb{R}^n }\{g(x^k,\epsilon^k)^Tx + \frac{\beta}{2}\|x - x^k\|^2\}, 
    \end{aligned}
    \end{equation}
    where each component of $g: \mathbb{R}^{2n}\to\mathbb{R}^n$ is defined by \eqref{smoothprob.grad}.
Therefore, $S$ is $C^1$ smooth in $B$. 

(iii) It follows from \eqref{sub.prob.eqv} that the iteration $T$ is $C^1$ smooth in $B$. 
\end{proof}

The above theorem shows the local behavior of ${\rm DIRL}_1$: it is equivalent to the gradient descent method for
the variables in $\Ical^*$    locally around a stationary point in $x^*+\mathcal{M}_{\Ical^*}$. Specifically, it can be expressed as 
     \begin{equation}\label{T.smooth}
     T_i(x,\epsilon) = \left\{\begin{aligned}
        &(1-\alpha)x_i, && i\in\Jcal^*, \\
        &x_i - \frac{\alpha}{\beta}\nabla_{x_i} F(x,\epsilon), && i\in\mathcal{I}^*,\\
        & (1-\alpha(1-\mu))\epsilon_i, && i \in \{n+1, \ldots, 2n\}
    \end{aligned}\right.
\end{equation} 
around $(x^*;\epsilon^*)$.  
We then proceed to the saddle escaping mechanism for DIRL$_1$, ensuring that the algorithm does not converge to the strict saddle points defined in Definition \ref{defSaddle}.  
To apply the CSM theorem, we have to show that the unstable fixed point of $T$ corresponds to
the strict saddle point of \eqref{pp}.

\begin{theorem}(Unstable Fixed Point) Suppose $x^*$ is a stationary point of \eqref{pp}.  
 $x^*$ is a strict saddle point of \eqref{pp} if and only if the corresponding $(x^*,\epsilon^*)$ with $\epsilon^*= 0$ is an unstable fixed point of $T$. 
 Moreover,  if $\alpha <  \frac{\beta}{\rho}$ where $\rho$ is defined by \eqref{hess is bounded},  $DT(x^*,\epsilon^*)$ is invertible. \label{ThLs2}
\end{theorem} 
\begin{proof} 
First of all, we calculate the Jacobian of $T$ in the neighborhood $B$ described in Theorem \ref{ModelId}. 
  It holds that 
        $D_x g(x^*,\epsilon^*) =
         \begin{bmatrix}  \nabla^2_{\Ical^*\Ical^*}  F(x^*) &  \nabla^2_{\Ical^* \Jcal^*}   f(x^*)   \\  \boldsymbol{0} & \beta I \end{bmatrix} $. 
   Therefore, the Jacobian of $T$ at $(x^*,\epsilon^*)$ is given by $$\begin{aligned}
    DT(x^*,\epsilon^*) = & \  \begin{bmatrix} (1-\alpha)I + \alpha(I-\tfrac{1}{\beta} D_x g(x^*,\epsilon^*))\  &D_{\epsilon} g(x^*,\epsilon^*) \\
        \boldsymbol{0} \ &   (1-\alpha(1-\mu)) I \end{bmatrix} \\
        = & \  \begin{bmatrix}
          I - \frac{\alpha}{\beta} \nabla^2_{\Ical^*\Ical^*}  F(x^*)  &  I - \frac{\alpha}{\beta}\nabla^2_{\Ical^* \Jcal^*}   f(x^*)   &  D_{\epsilon} g(x^*,\epsilon^*)  \\ 
              \boldsymbol{0}    &  (1-\alpha) I  &   \boldsymbol{0}   \\ 
              \boldsymbol{0}  &\boldsymbol{0} &     (1-\alpha(1-\mu)) I 
        \end{bmatrix} 
        \end{aligned}$$

 Now we are ready to prove the equivalence between the unstable fixed point of $T $ and the strict saddle point of \eqref{pp}.      Recall Definition \ref{def.fix.point} and Definition \ref{defSaddle} and notice $\mu \in (0,1)$ and $\alpha\in(0,1)$.   It suffices to show that $I - \frac{\alpha}{\beta} \nabla^2_{\Ical^*\Ical^*}  F(x^*)  $ has an eigenvalue greater than 1 if and only if $\nabla^2_{\mathcal{I}^*\Ical^*}  F(x^*)$ has a negative eigenvalue.  
 This equivalence is obviously true since 
     \[\lambda_{\min}( \nabla^2_{\mathcal{I}^*\Ical^*}  F(x^*) ) < 0 
     \iff  \lambda_{\max}(  I - \frac{\alpha}{\beta} \nabla^2_{\Ical^*\Ical^*}  F(x^*)  ) > 1.\]
     
Moreover, if $\alpha$ is selected to satisfy $1-\frac{\alpha}{\beta}\rho > 0$ and by the strict saddle property 
$ \nabla^2_{\mathcal{I}^*\Ical^*}  F(x^*)$ does not have a 0 eigenvalue,  we know by \eqref{hess is bounded} that 
$DT(x^*,\epsilon^*)$ is invertible. 
 \end{proof}

To apply Corollary \ref{coro.saddle}, we next show that $T$ is a lipeomorphism, even if it may not be differentiable globally in $\mathbb{R}^{2n}$. This is equivalent to showing that $T$ and its inverse $T^{-1}$ are Lipschitz continuous.

\begin{proposition}(Lipeomorphism) \label{LmaLip} Suppose Assumptions \ref{ass.1}-\ref{ass.3} hold.  
     Let $L_r :=  r''(| (r')^{-1}(C)|)$ if $r'(|0^+|) = \infty$, where $C$ is defined in \eqref{xnablabd},    and $L_r :=  r''(|0^+|)$ if $r'(|0^+|) < \infty$.   
  Then $T$ is invertible; moreover, $T$ and $T^{-1}$ are Lipschitz continuous (that is, $T$ is a lipeomorphism) when $0 < \alpha < 1$,  $\beta > \alpha L_{\nabla f}/2$ and  $\alpha(2  +  \tfrac{1}{\beta} L_{\nabla f} + \tfrac{\lambda}{\beta} L_r +\mu)  < 1$.  
\end{proposition}
\begin{proof} 


If $r'(|0^+|) < \infty$,  $w(x_i,\epsilon_i) \le  r'(|0^+|)$ by Assumption \ref{ass.1} and $w(x_i,\epsilon_i) =  r'(|x_i|+\epsilon_i)$ is  Lipschitz continuous  for $x\in \Lcal(F^0)$ and $ \epsilon_i\in\mathbb{R}_+$ by Assumption \ref{ass.2}. 

If $r'(|0^+|) = \infty$,  by \eqref{xnablabd},
$  \| x- \tfrac{1}{\beta} \nabla f(x)\|_{\infty} \leq C$, $\forall x\in \Lcal(F^0)$. 
 When $  w(x_i,\epsilon_i) > C/\lambda$,  by the definition of $S^x$, 
 \[  0 \in   \nabla_i f(x) + \beta (y_i -  x_i ) + \lambda  w(x_i,\epsilon_i) \partial |y_i|. \]
It is true that $S_i^x(x,\epsilon)= 0$. Therefore,  for any $x\in \Lcal(F^0)$ with $S_i^x(x,\epsilon)\ne 0$,  
\begin{equation}\label{eq.bound} 
 r'(|x_i|+\epsilon_i) \le C/\lambda \text{ \ and \ }  |x_i|+\epsilon_i \ge \underline{x} := (r')^{-1}(C/\lambda) > 0
 \end{equation}
  by Assumption \ref{ass.1}.

It  implies that 
 we can equivalently write $S$ as 
 \[ S_i^x(x,\epsilon) = \mathbb{S}_{\lambda  u_i   (x_i,\epsilon_i)  /\beta} ( x_i - \tfrac{1}{\beta}\nabla_i f(x)),\]
 where we define  $u: \mathbb{R}^{2n}\to \mathbb{R}^n$ with 
 \[  u_i (x_i,\epsilon_i) =
 \begin{cases} 
  w(x_i,\epsilon_i)= r'(|x_i|+\epsilon_i) & \text{  if   }   r'(|0^+|) < \infty,\\
   \min( C/\lambda, w(x_i,\epsilon_i)) 
  = \min(C/\lambda,  
  r'(|x_i|+\epsilon_i)) & \text{  if   }   r'(|0^+|) = \infty.
 \end{cases}\] 
 Therefore,   $u(\cdot)$ is also Lipschitz continuous for  $x\in\Lcal(F^0)$ and $   \epsilon \in \mathbb{R}^n_+$.  Let $L_r$ be its Lipschitz constant, meaning for any $\hat x, \tilde x \in  \Lcal(F^0)$ and $\hat \epsilon, \tilde\epsilon \in\mathbb{R}^n_+$ 
  \begin{equation}\label{u.nonexp} \| u(\hat x, \hat \epsilon) - u(\tilde x, \tilde \epsilon)  \| \le L_r \|  (\hat x; \hat \epsilon) - (\tilde x; \tilde \epsilon) \|.\end{equation}
  A natural  choice   is 
 $L_r= r''(|0^+|)$ if $r'(|0^+|)<\infty$ and $L_r = r''(|\underline{x}|)$ if $r'(|0^+|)=\infty$.

 Next, we prove that $S$ is Lipschitz continuous.  It follows that 
 \[ \begin{aligned} 
 {\|} S^x(\hat x, \hat \epsilon) - S^x(\tilde x, \tilde \epsilon) {\|} 
  = &\   {\|} \Sb_{\lambda \hat u/\beta}(\hat x - \tfrac{1}{\beta}\nabla f(\hat x)) - \Sb_{\lambda \tilde u/\beta}(\tilde x - \tfrac{1}{\beta}\nabla f(\tilde x)) {\|} \\ 
   \le  &\     \| \hat x - \tfrac{1}{\beta}\nabla f(\hat x)- (\tilde x - \tfrac{1}{\beta}\nabla f(\tilde x))\| +  \tfrac{\lambda}{\beta} \| \hat u - \tilde u \|  \\
   \le & \     \| \hat x - \tilde x \| +  \tfrac{1}{\beta} \|\nabla f(\hat x) -  \nabla f(\tilde x))\| + \tfrac{\lambda}{\beta} L_r \|  (\hat x; \hat \epsilon) - (\tilde x; \tilde \epsilon) \|  \\
   \le & \ (1 +  \tfrac{1}{\beta} L_{\nabla f}) \| \hat x  -   \tilde x \| + \tfrac{\lambda}{\beta} L_r \|  (\hat x; \hat \epsilon) - (\tilde x; \tilde \epsilon) \|\\
   \le & \ (1 +  \tfrac{1}{\beta} L_{\nabla f} + \tfrac{\lambda}{\beta} L_r ) \|  (\hat x; \hat \epsilon) - (\tilde x; \tilde \epsilon) \|, 
   \end{aligned} 
   \]
  where the first inequality is determined by  Lemma \ref{lem.nonexpansive} and the second is by \eqref{f.nonexp} and \eqref{u.nonexp}. 

\noindent Therefore, for any $\hat x, \tilde x \in \Lcal(F^0)$ and $\hat\epsilon,\tilde\epsilon \in\mathbb{R}^n_{+}$, 
\[ \begin{aligned}
  & \| \alpha(  I- S) (\hat x,\hat \epsilon)     -  \alpha(  I- S) (\tilde x,\tilde \epsilon)   \| \\
  \le & \ \alpha \| (\hat x;\hat \epsilon) -(\tilde x;\tilde \epsilon)\| + \alpha \| S^x (\hat x,\hat \epsilon) -  S^x (\tilde x,\tilde \epsilon)  \| + \alpha \| S^\epsilon (\hat x,\hat \epsilon) -  S^\epsilon (\tilde x,\tilde \epsilon)  \|  \\
  \le &\ \alpha(2  +  \tfrac{1}{\beta} L_{\nabla f} + \tfrac{\lambda}{\beta} L_r ) \|  (\hat x; \hat \epsilon) - (\tilde x; \tilde \epsilon) \|+
  \alpha\mu \| \hat \epsilon - \tilde \epsilon\|\\
    \le &\ \alpha(2  +  \tfrac{1}{\beta} L_{\nabla f} + \tfrac{\lambda}{\beta} L_r +\mu) \|  (\hat x; \hat \epsilon) - (\tilde x; \tilde \epsilon) \|. 
  \end{aligned} \]


\noindent Notice that 
     $T =  I + \alpha(S - I)$. 
    Since $  \alpha(I-S) $ is Lipschitz with constant $\alpha(2  +  \tfrac{1}{\beta} L_{\nabla f} + \tfrac{\lambda}{\beta} L_r +\mu)  < 1$.  Applying  \cite[Lemma 2.13]{davis2022proximal}, we directly derive that  $T = I + \alpha(S - I)$ is invertible and a lipeomorphism.
\end{proof}

Now we are ready to apply Corollary \ref{coro.saddle} to derive the following global escape theorem.
\begin{theorem}(Global Escape)
    Consider problem \eqref{pp} with strict saddle property and the ${\rm DIRL}_1$ algorithm with Assumptions \ref{ass.1}-\ref{ass.3} being satisfied. Suppose  $0 < \alpha < 1$,  $\beta > \alpha L_{\nabla f}/2$, $\alpha(2  +  \tfrac{1}{\beta} L_{\nabla f} + \tfrac{\lambda}{\beta} L_r +\mu)  < 1$ and $\alpha < \beta/\rho$. Let $\{(x^k,\epsilon^k)\}$ be the sequence generated by the Algorithm \ref{alg2}. Then, the set of initial points attracted by strict saddle points has zero Lebesgue measure. 
\end{theorem}
\begin{proof}
    It follows directly from Theorem \ref{ThLs2} that the map $T$ is local $C^1$ smooth around any stationary point, and that all stationary points are unstable fixed points. 
    By applying Proposition \ref{LmaLip}, $T$ is a lipeomorphism. 
    By Corollary \ref{coro.saddle},   the set of initial conditions attracted by such fixed points has a zero Lebesgue measure. 
\end{proof}

\section{Damped Iterative Reweighted $\ell_2$ Algorithms} \label{sec4}
The ${\rm IRL}_2$ algorithm is an effective method to solve problem \eqref{pp}.  In this method, we write
$h(\cdot ) = r\circ \sqrt{(\cdot)}$, so that $s(x_i) = h(x^2_i)$.  
Since 
 $r$ is increasing and concave, it follows that 
$h$ is also concave on $\mathbb{R}_+$.

At each iteration, after adding perturbation $\epsilon_i^2$ to each $h$ and linearizing $h$ at $(x_i^k)^2 + (\epsilon_i^k)^2$,  the ${\rm IRL}_2$ algorithm solves the following weighted $\ell_2$ regularized subproblem (for brevity, we still use 
$S$ and $T$ for the maps),  
\begin{align*}
 y^k =   S^x(x^k,\epsilon^k)  = \arg\min_{y} G(y;x^k):=\nabla f(x^k)^Ty + \tfrac{\beta}{2}\|x^k - y\|^2 + \lambda\sum_{i=1}^n u(x^k_i,\epsilon_i^k) y_i^2,
\end{align*}
where 
\[ u(x_i,\epsilon_i) = h'( x_i^2+\epsilon_i^2 ) =  \frac{r'(\sqrt{x_i^2 + \epsilon_i^2})}{2\sqrt{x_i^2 + \epsilon_i^2}} \ \text{  and } \ 
 u(x_i, 0)  = \begin{cases}  \frac{r'(   |x_i|  )}{2 |x_i|  }   & \text{ if }  x_i \ne 0 \\ 
  \infty & \text{ if }  x_i=0.  
 \end{cases}
 \]
 The solution of the subproblem is given by   
 \begin{equation}
   y^k_i =   S_i^x(x,\epsilon) := \frac{1}{1 + \frac{2\lambda}{\beta}u(x_i^k,\epsilon_i^k)}(x_i^k - \tfrac{1}{\beta}\nabla_i f(x^k)).
\end{equation}
    If $u(x_i^k,\epsilon_i^k) = \infty$, we simply set $S_i^x(x^k,\epsilon^k) = 0$. The damped iteratively reweighted $\ell_2$ algorithm, hereinafter named DIRL$_2$ algorithm,  is presented in Algorithm \ref{alg2}.  
 \begin{algorithm}
    \begin{algorithmic}[1]
        \caption{Damped Iteratively Reweighted $\ell_2$ Algorithm (${\rm DIRL}_2$)}
        \label{alg2}
        \STATE \textbf{Input} $\beta > \alpha L_{\nabla f}/2$, $\alpha\in(0, 1)$, $x_0$, $\epsilon_0 > 0$, and $\mu \in (0,1)$.
        \STATE \textbf{Initialize}: set $k = 0$.
        \REPEAT 
            \STATE \textrm{Reweighting:} $u(x_i^k, \epsilon_i^k) = h'((x_i^k)^2 + (\epsilon_i^k)^2)$.
            \STATE \textrm{Compute new iterate:} 
            \STATE $y^{k} = \arg\min\limits_{y\in\mathbb{R}^n} G(y;x^k):=\nabla f(x^k)^Ty + \tfrac{\beta}{2}\|x^k - y\|^2 + \lambda\sum_{i=1}^n u(x^k_i,\epsilon_i^k) y_i^2.$
            \STATE $x^{k+1} = (1-\alpha)x^k + \alpha y^{k}$.
            \STATE $\epsilon^{k+1} = \mu\epsilon^k$. 
            \STATE \textrm{Set} $k = k + 1$.
        \UNTIL{Convergence}
    \end{algorithmic}
\end{algorithm} 

The fixed point iteration of DIRL$_2$ is then given as 
$(x^{k+1}; \epsilon^{k+1}) = T(x^k, \epsilon^k)$ with $T = (1-\alpha)I + \alpha S$ and 
\begin{equation}
      T (x,\epsilon) := \begin{bmatrix}
        (1-\alpha)x + \alpha S^x(x,\epsilon) \\ 
          (1-\alpha(1-\mu)) \epsilon
    \end{bmatrix}.
\end{equation}

We have the following result to connect  the stationary point of \eqref{pp} with the fixed point of $T$. 
\begin{proposition}\label{prop.fix.point.2} 
 If $x^*$ is a stationary point of \eqref{pp}, then  $x^* = S^x(x^*,0)$  and  $(x^*; 0) = T(x^*, 0)$. Conversely, 
 for any $x^*$ satisfying $x^* = S^x(x^*,0)$,  
 if $r'(|0^+|)=\infty$ or $\Jcal^*=\emptyset$,  then $x^*$ is a stationary point of \eqref{pp}. 
\end{proposition} 
\begin{proof} 


Suppose $x^*$ is a stationary point. 
It follows from \eqref{opt.cond1} that  for $i\in\Ical^*$,  
$x_i^*  = x_i^* - \tfrac{1}{\beta}( \nabla_if(x^*) + \lambda r'(|x_i^*|)\sign(x_i^*) ) $ implying 
$ x_i^*\left( 1 + \tfrac{2\lambda}{\beta}  \tfrac{r'(|x_i^*|)}{2|x_i^*|} \right) = x_i^* - \tfrac{1}{\beta} \nabla_if(x^*)$. 
Therefore,  $ x_i^* =  S_i^x(x^*, 0 )$ since $u(x_i^*, 0) = \tfrac{r'(|x_i^*|)}{2|x_i^*|}$.  
As for $i \in \Jcal^*$, $w(0,0) = \infty $, implying that 
$ 0 = S_i(x^*, 0)$. In summary, $x^*= S^x(x^*, 0)$.

Now consider $x^* = S^x(x^*,0)$. If $r'(|0^+|)=\infty$ or $\Jcal^*=\emptyset$, it suffices to verify 
\eqref{opt.cond1} is satisfied by $x^*$.  
It follows from $x^*_i = S_i^x(x^*,0), i\in \Ical^*$ that 
$ x_i^*\left( 1 + \tfrac{2\lambda}{\beta}  \tfrac{r'(|x_i^*|)}{2|x_i^*|} \right) = x_i^* - \tfrac{1}{\beta} \nabla_if(x^*)$, implying 
$\nabla_if(x^*) + \lambda r'(|x_i^*|)\sign(x_i^*) = 0$, completing the proof. 
\end{proof}

\begin{Remark}
The converse statement of Proposition \ref{prop.fix.point.2} may not be   true when $\Jcal^*\ne \emptyset$ 
and $r'(|0^+|) < \infty$. Since for any $x^*$ that satisfies
$x_i^* = S_i^x(x^*,0)$, $i\in\Ical^*$, it is naturally true that $u(0,0) =\infty$ and $S_i^x(x^*,0) = 0$. However, there is no guarantee that $x^*$ satisfies \eqref{opt.cond2}.
    In those cases,  $S$ may have two types of fixed points --- those are stationary points and others are not. 
 This represents a stark difference from the DIRL$_1$ algorithm, where all fixed points of the iteration map are stationary points. 
\end{Remark}

%
%
%

\subsection{Convergence analysis}

In this subsection, we show that DIRL$_2$ only converges to its fixed points which are the stationary points of \eqref{pp} under the 
assumption that   $r'(|0^+|)=\infty$ or $\Jcal^*=\emptyset$. 
 Similarly to the analysis of the DIRL$_1$ algorithm, the analysis of DIRL$_2$ is also based on the monotonic decrease in the objective
$F(x, \epsilon)$.    We also need Assumption \ref{ass.3} and consequently \eqref{f.nonexp} to hold for the analysis. However, 
it should be noted that the definition of $F(x,\epsilon)$ is different from that of the DIRL$_1$ algorithm.

\begin{proposition} 
    Suppose Assumptions \ref{ass.1}-\ref{ass.3} hold. Let $\{(x^k, \epsilon^k)\}$ be the sequence generated by Algorithm \ref{alg2} with $\alpha \in (0,1)$ and   $\frac{\beta}{\alpha} \geq \frac{L_{\nabla f}}{2}$. It holds that \{$F(x^k, \epsilon^k)$\} is monotonically decreasing and the reduction satisfies
    \begin{equation*}
        F(x^0, \epsilon^0) - F(x^k, \epsilon^k) \geq (\frac{\beta}{\alpha} - \frac{L_{\nabla f}}{2})\sum\limits_{t=0}^{k-1}\|x^t - x^{t+1}\|^2.
    \end{equation*}
  Moreover, $ \lim\limits_{k\rightarrow\infty} \|x^{k+1} - x^{k}\| = 0$ and $\{x^k\} \subset \Lcal(F^0)$. 
\end{proposition}
\begin{proof} 
    It follows from the concavity of $h$ that   for $i = 1, ..., n$,
         \begin{equation*}
        h(x_i^{k+1})^2 +(\epsilon_i^k)^2 ) \leq h(x_i^k)^2+(\epsilon_i^k)^2 ) + h'(( x_i^k)^2+(\epsilon^k_i)^2)((x_i^{k+1})^2 - (x_i^k)^2), 
    \end{equation*}
which is equivalent to 
    \begin{equation*}
        r( \sqrt{ (x_i^{k+1})^2 +(\epsilon_i^k)^2}) \leq r(\sqrt{ (x_i^k)^2+(\epsilon_i^k)^2}) + u( x_i^k, \epsilon^k_i)((x_i^{k+1})^2 - (x_i^k)^2).
    \end{equation*}
    Following the same argument for \eqref{p2.con2}-\eqref{p2.con4}, we have 
\begin{equation}
   \Delta F(x^{k+1}, \epsilon^{k+1}) \geq  \Delta G(x^{k+1}; x^k, \epsilon^k) + \frac{\beta - L_{\nabla f}}{2}\|x^k - x^{k+1}\|^2. \label{p2.con4.l2}
  \end{equation}
     
    For each subproblem, $y^k$ satisfies the optimality condition
    \begin{equation}
        \nabla f(x^k) + 2 \lambda U(x^k, \epsilon^k)y^k + \beta(y^k - x^k)= 0.  \label{p2.opt.l2}
    \end{equation}
     where $U(x^k, \epsilon^k) := {\rm diag}(u(x_1^k, \epsilon_1^k), ..., u(x_n^k, \epsilon_n^k))$. Consequently, by $x^{k+1} = (1-\alpha)x^k + \alpha y^k$ and the optimality condition \eqref{p2.opt}, the following holds 
    \begin{equation} \label{p2.con5}
    \begin{aligned}
     & \ \Delta G(x^{k+1}; x^k, \epsilon^k) \\
    =&\ \nabla f(x^k)^T(x^k - x^{k+1}) + \lambda\sum\limits_{i=1}^{n}u(x_i^k,\epsilon_i^k)((x_i^k)^2 - (x_i^{k+1})^2) - \frac{\beta}{2}\|x^k - x^{k+1}\|^2 \\
        \ge &\ \nabla f(x^k)^T(\alpha x^k - \alpha y^k) + \lambda\sum\limits_{i=1}^{n}u(x_i^k,\epsilon_i^k)(\alpha(x_i^k)^2 - \alpha(y_i^k)^2) - \frac{\beta\alpha^2}{2}\|x^k - y^k\|^2 \\
        \geq&\ \alpha [ \nabla f(x^k)^T(x^k - y^k) + 2\lambda\sum\limits_{i=1}^{n}u(x_i^k, \epsilon_i^k)y_i^{k}(x_i^k - y_i^k)] - \frac{\beta\alpha^2}{2}\|x^k - y^k\|^2 \\
        =&\ \alpha[\nabla f(y^k) + 2 \lambda U(x^k, \epsilon^k)y^{k} + \beta(y^k - x^k)]^T(x^k - y^k)  + \frac{\alpha\beta(2 - \alpha)}{2}\|x^k - y^k\|^2 \\
        =&\ (\frac{\beta}{\alpha} - \frac{\beta}{2})\|x^k - x^{k+1}\|^2. 
    \end{aligned}
    \end{equation} 
            Here the first inequality is by $(1-\alpha)(x_i^k)^2 - \alpha  ( y_i^k)^2 \ge (\alpha x_i^k + (1-\alpha)y_i^k )^2$ 
             and the second is by the convexity of $(\cdot)^2$.  
    Combining \eqref{p2.con4} and \eqref{p2.con5}, we get
    \begin{equation*}
        \Delta F(x^{k+1}, \epsilon^{k+1}) \geq (\frac{\beta}{\alpha} - \frac{L_{\nabla f}}{2} )\|x^k - x^{k+1}\|^2.
    \end{equation*}
    Replacing $k$ with $t$ and summing up from $t = 0$ to $k-1$, we have
    \begin{equation*}
        \sum\limits_{t=0}^{k-1}[F(x^t, \epsilon^t) - F(x^{t+1}, \epsilon^{t+1})] = F(x^0, \epsilon^0) - F(x^k, \epsilon^k) \geq (\frac{\beta}{\alpha} - \frac{L_{\nabla f}}{2})\sum\limits_{t=0}^{k-1}\|x^t - x^{t+1}\|^2.
    \end{equation*}
    It follows that $F(x^k) \le F(x^k,\epsilon^k) \le F^0$, which implies that $\{x^k\} \in \mathcal{L}(F^0)$. 
    By Assumption \ref{ass.3}, it is bounded 
    and $\mathop{\lim\inf}\limits_{k \rightarrow \infty} F(x^k, \epsilon^k) > -\infty$. Taking $k \rightarrow \infty$, we have  
    \begin{equation*}
        (\frac{\beta}{\alpha} - \frac{L_{\nabla f}}{2})\sum\limits_{t=0}^{\infty}\|x^t - x^{t+1}\|^2 \leq F(x^0, \epsilon^0) - \mathop{\lim\inf}\limits_{k \rightarrow \infty} F(x^k, \epsilon^k) < \infty.
    \end{equation*}
    which implies $\lim\limits_{k \rightarrow \infty}\|x^k - x^{k+1}\| = 0$.
\end{proof}

With the monotonicity of $F$, we derive the global convergence for the ${\rm DIRL}_2$ method. 
\begin{theorem}(Global convergence)\label{GlobalConverge.p2}
    Suppose that Assumptions \ref{ass.1}--\ref{ass.3} hold, and $\Jcal^*=\emptyset$ or $r'(|0^+|)=\infty$. Let $\{(x^k,\epsilon^k)\}$ being generated by ${\rm DIRL}_2$ and $\Omega$ be the set of the limit points of $\{x^k\}$. Then any $x^*\in \Omega$ is a stationary point for \eqref{pp}.  
\end{theorem}

\begin{proof}
Firstly, the boundedness of $\mathcal{L}(F^0)$ from Assumption \ref{ass.3} implies $\Omega \neq \emptyset$. Then let $x^*$ be a limit point of $\{x^k\}$ with subsequence $\{x^k\}_{\mathcal{S}} \rightarrow x^*$ and $\{y^k\}_{\Scal} \rightarrow x^*$ since $y^k = \frac{1}{\alpha}(x^k - x^{k+1}) + x^k$.   
The optimality condition of the $k$th subproblem implies 
\begin{align}\label{l2.res.0}
 - \beta(y_i^k - x_i^k) = \nabla_if(x^{k})  + 2 \lambda u(x_i^k, \epsilon_i^k)   y_i^k = \nabla_if(x^k) +  \frac{r'(\sqrt{(x_i^k)^2 + (\epsilon_i^k)^2})}{\sqrt{(x_i^k)^2 + (\epsilon_i^k)^2}} y_i^k.
\end{align} 
  As $  x_i^k \overset{\Scal}{\to} x_i^*$, $  y_i^k \overset{\Scal}{\to} x_i^*$ and 
 $\epsilon_i^k \to 0$,   $ \lim\limits_{k\in \Scal}\frac{y_i^k}{\sqrt{x_i^2 + \epsilon_i^2}} = \sign(x_i^*)$ for $i\in \Ical^*$. 
 Therefore, \eqref{opt.cond1} is satisfied. 

If $\Jcal^* = \emptyset $, there is nothing to prove. If not, for $i\in\Jcal^*$, there is nothing to prove  since $|\nabla_i f(x^*)| \le r'(|0^+|)$ naturally holds. 
     In general, $x^*$ is a stationary point.
  \end{proof}

 \subsection{Avoidance of strictly saddle points} 
 
 In this subsection, we further show that the DIRL$_2$ algorithm does not converge to fixed points that are stationary points and strict saddle points. Unlike DIRL$_1$,  the map $S$ in DIRL$_2$ does not enjoy the property of ``active manifold identification'' or ``model identification''.  
There is no guarantee that the iteration map is $C^1$-smooth without this property. The following additional assumption is needed to 
guarantee local smoothness, which is satisfied by the $\ell_p$ regularization. 

\begin{assumption} \label{ass.4} 
$\lim\limits_{z\rightarrow 0^+} r'(z)  =\infty\   \text{ and } \        \lim\limits_{z\rightarrow 0^+}z\frac{ r''(z)  }{r'(z)^2} \to 0.$ 
\end{assumption} 

 To show the smoothness of $S$, we first calculate the Jacobian of $S$.  
Define the auxiliary  function  $g: \mathbb{R}_+\to \mathbb{R}$: 
\[ g(z) : =\frac{1}{1+\tfrac{\lambda}{\beta}\tfrac{r'(z)}{z}} \ \text{ with } g(0) := 0.\] 
We can write  
$S_i^x(x,\epsilon) = g(z)(x_i - \tfrac{1}{\beta}\nabla_i f(x))$ with 
  $z =\sqrt{x^2 + \epsilon^2}$. 
  It follows that 
\begin{equation} \label{eq.g.grad} 
g'(z) = 
  \begin{cases}  -\tfrac{\lambda}{\beta}\frac{r''(z)z -  r'(z)}{(z +\tfrac{\lambda}{\beta} {r'(z)} )^2}, & \text{\  for }  z\ne 0, \\
  \lim\limits_{\delta\to 0} \left(\frac{1}{1+\tfrac{\lambda}{\beta}\tfrac{r'(|\delta|)}{|\delta |}} - 0 \right) \frac{1}{\delta}=  0, &  \text{\  for }  z= 0. 
  \end{cases}\end{equation}
   It holds that 
 \[ \begin{aligned} \lim\limits_{z\to 0}  g'(z) = \ &
-\tfrac{\lambda}{\beta} \lim\limits_{z\to 0} \left(  \frac{r''(z)z}{(z + \tfrac{\lambda}{\beta} {r'(z)} )^2} - \frac{  r'(z)}{(z + \tfrac{\lambda}{\beta} {r'(z)} )^2}   \right)\\ 
  = \ &- \lim\limits_{z\to 0}   \left(   \tfrac{\beta}{\lambda} \frac{r''(z)z}{  {r'(z)} ^2} -   \frac{  z  }{   r'(z)    } \right)
  = 0 = g'(0)
  \end{aligned}\]
  by Assumption \ref{ass.4}. Therefore, $g$ is continuously differentiable in $\mathbb{R}$. 


 We have the following result for maps $S$ and $T$. 
 
\begin{theorem}\label{ModelId2}    
 (Local smoothness and Lipschitz continuity) Suppose Assumptions \ref{ass.1}-\ref{ass.4} hold and  $x^*$ is a stationary point of \eqref{pp}. 
 Then the following statements hold: 
 \begin{enumerate}
  \item[(i)] $S$ is $C^1$ smooth and Lipschitz continuous on $\Lcal(F^0)$ with constant $L_S$; moreover  
  \begin{equation}\label{jac.h} 
 DS(x^*, \boldsymbol{0})=\   \begin{bmatrix}    H_{\Ical^* \Ical ^*}  &  H_{\Ical^* \Jcal^* } & \boldsymbol{0} \\ 
\boldsymbol{0}&\boldsymbol{0} & \boldsymbol{0} \\ 
\boldsymbol{0}&\boldsymbol{0} & \mu I \end{bmatrix},    
 \end{equation}
where 
 $H_{\Ical^* \Ical ^*} = I-  \tfrac{1}{\beta} (I+ \tfrac{\lambda}{\beta}U_{\Ical^* \Ical ^*})^{-1} \nabla^2_{\Ical^* \Ical ^*} F(x^* )$.  
  \item[(ii)] $T$ is a lipeomorphism for $\alpha \in(0, \frac{1}{1+L_S})$ with 
 \[  DT(x^* ,  \boldsymbol{0})=\   \begin{bmatrix}   (1-\alpha)I  + \alpha H_{\Ical^* \Ical^* }  & \alpha H_{\Ical^* \Jcal^* } &  \boldsymbol{0} \\ 
  \boldsymbol{0}&(1-\alpha)I &   \boldsymbol{0}  \\ 
  \boldsymbol{0} &  \boldsymbol{0}   &  (1-\alpha(1-\mu )) I \end{bmatrix}.\]
\end{enumerate} 
\end{theorem}

   \begin{proof}   
  (i)  We simply calculate its partial derivative and verify that it is continuous at any point $(x^*; \epsilon^*) \in \Rb^{2n}$.  
  The smoothness of $S^\epsilon$ is obvious. We only verify the smoothness of $S^x$.  

Consider $i$ satisfying $z_i^* : = \sqrt{(x_i^*)^2 + (\epsilon_i^*)^2} > 0$. In this case,  $S$ is smooth with respect to 
 $x_i$ and $\epsilon_i$, and  
   \begin{align} 
\frac{\partial}{\partial x_i} S_i^x(x^*, \epsilon^*)  = & \  g(z^*_i)(1-\tfrac{1}{\beta}\nabla^2_{ii} f(x^*)) + \tfrac{x_i^*}{z_i^*} g'(z_i^*)(x_i^*-\tfrac{1}{\beta}\nabla_i f(x^*)); \label{jac.ii}\\
\frac{\partial}{\partial x_j} S_i^x(x^*, \epsilon^*)  = & \    
- g(z^*_i) \tfrac{1}{\beta}\nabla^2_{ij} f (x^*),  \quad  j\ne i; \label{jac.ij}\\
\frac{\partial}{\partial \epsilon_i} S_i^x(x^*, \epsilon^*)  = & \      \tfrac{ \epsilon_i^* }{z_i^*} g'(z^*_i)(x_i^* - \tfrac{1}{\beta} \nabla_i f (x^*)); \label{jac.ieps}\\
\frac{\partial}{\partial \epsilon_j} S_i^x(x^*, \epsilon^*)  = &\  0, \quad  j\ne i.\label{jac.jeps}
\end{align}

Consider $ i  $ satisfying $ x_i^*=0$ and $\epsilon_i^*   = 0$.  It follows from $S_i^x(x ^*, \boldsymbol{0})=0$ and  $r'(|0^+|) = \lim\limits_{\delta \to 0} u(  \delta, 0) = \lim\limits_{\delta \to 0} u( 0, \delta)= \infty$  that 
\begin{equation}\label{jac.ii1}
     \begin{aligned}
    \frac{\partial  }{\partial x_i} S _i^x(x^* , \epsilon^* )  = & \ \lim\limits_{\delta\rightarrow 0}  \frac{1}{1 + \frac{2\lambda}{\beta}u(  \delta, 0)} [  \delta - \tfrac{1}{\beta}\nabla_i f(x^*  + \delta e_i )]/\delta \\ 
        =&\   - \tfrac{1}{\beta}\nabla_i f(x^* )\lim\limits_{\delta\rightarrow 0}[g(\delta)-g(0)] /\delta  \\
                =&- \tfrac{1}{\beta}\nabla_i f(x^*  )\lim\limits_{\delta\rightarrow 0}g'(0)            = 0, \end{aligned}
\end{equation}

        \begin{equation}\label{jac.ii2}
        \begin{aligned}
  \frac{\partial  }{\partial x_j}S^x_i (x^*,\epsilon^*)  = & \ \lim\limits_{\delta\rightarrow 0}  \frac{1}{1 + \frac{2\lambda}{\beta}u(  \delta, 0)} [  \delta - \tfrac{1}{\beta}\nabla_i f(x^*  + \delta e_j )]/\delta  \\ 
          =&\   - \tfrac{1}{\beta}\nabla_i f(x^* )\lim\limits_{\delta\rightarrow 0}[g(\delta)-g(0)] /\delta     =     0, \quad   j\ne i; 
              \end{aligned}\end{equation}
              
\begin{equation}\label{jac.ii3}
                \begin{aligned} 
  \frac{\partial  }{\partial \epsilon_i} S _i^x(x^* , \epsilon^* )  = & \ \lim\limits_{\delta\rightarrow 0}  \frac{1}{1 + \frac{2\lambda}{\beta}u( 0, \delta )} [    - \tfrac{1}{\beta}\nabla_i f(x ^*)]/\delta  \\
        =&\   - \tfrac{1}{\beta}\nabla_i f(x^* )\lim\limits_{\delta\rightarrow 0}[g(\delta)-g(0)] /\delta     =     0;
        \end{aligned}\end{equation}
              
                  \begin{equation}\label{jac.ii4}
   \frac{\partial  }{\partial \epsilon_j}S^x_i (x^*,\epsilon^*)  =  \  0, \quad  j\ne i.  
 \end{equation}
 
\noindent Now we verify that the Jacobian is continuous at $(x^*, \epsilon^*)$. It suffices to show this  for any $ i  $ satisfying $ x_i^*=0$ and $\epsilon_i^*   = 0$. 
It follows from \eqref{jac.ii},  the fact that $\frac{x_i }{z_i } \le 1$ and $g(0)=0, g'(0)=0$ that 
 \[ \lim_{x\to x^*,\epsilon\to\epsilon^*} \frac{\partial}{\partial x_i} S_i^x(x, \epsilon) = 0 =  \frac{\partial}{\partial x_i} S_i^x(0 , 0 ).\]
It follows from \eqref{jac.ij} and    $g(0)=0$ that 
\[ \lim_{x\to x^*,\epsilon\to\epsilon^*}  \frac{\partial}{\partial x_j} S_i^x(x, \epsilon) = 0 =  \frac{\partial}{\partial x_j} S_i^x(0 , 0 ).\] 
It follows from \eqref{jac.ieps}, the fact that $\frac{x_i }{z_i } \le 1$ and $ g'(0)=0$ that  
  \[\lim_{x\to x^*,\epsilon\to\epsilon^*}  \frac{\partial}{\partial \epsilon_i} S_i^\epsilon(x, \epsilon)  =  0 =  \frac{\partial}{\partial \epsilon_i} S_i^\epsilon(0 , 0 ).\]
  It  follows from  \eqref{jac.jeps} that 
\[\lim_{x\to x^*,\epsilon\to\epsilon^*}  \frac{\partial}{\partial \epsilon_j} S_i^\epsilon(x, \epsilon)  =  0 =  \frac{\partial}{\partial \epsilon_j} S_i^\epsilon(0 , 0 ). \]
Overall, we have shown that $S$ is continuously differentiable.

To derive   \eqref{jac.h},  now suppose $x^*$ is a stationary point for \eqref{pp} and $\epsilon^*= \boldsymbol{0}$. 

For $i\in \Ical ^*$,  \eqref{opt.cond1} is equivalent to $ x_i^* - \frac{1}{\beta} \nabla_i f(x^*) = x_i^* + \frac{\lambda}{\beta} \frac{x_i^*}{z_i^*} r'(|x^*_i|)$ and $x_i^*/z_i^* = \text{sign}(x_i^*)$. It  follows from \eqref{jac.ii} that 
for $i\in \Ical ^*$, 
 \begin{equation}\label{eq.hii} 
 \begin{aligned}
     h_{ii}   =  & \   \frac{\partial}{\partial x_j} S^x_i(x ^* ,\boldsymbol{ 0})  \\
=   & \    \frac{(1 -  \tfrac{1}{\beta}\nabla^2_{ii} f(x ^* ))  }{  1 +  \tfrac{\lambda}{\beta}\frac{r'(|x ^* _i|)}{|x ^* _i|}     } - \frac{   (|x ^*_i | + \tfrac{\lambda}{\beta} r'(x  ^*_i)) \tfrac{\lambda}{\beta} \left(  |x  ^*_i| r''(|x  ^*_i|) -  r'(|x ^* _i|) \right)}{\left( |x  ^*_i| +  \tfrac{\lambda}{\beta} r'(|x ^* _i|)  \right)^2},\\
      =   & \   1- \tfrac{1}{\beta} \frac{  \nabla^2_{ii} f(x  ^*)   +  \lambda  r''(|x _i ^*|)   }{ 1+  \tfrac{\lambda}{\beta} \frac{r'(|x ^* _i|)}{|x  ^*_i|}}. \end{aligned}\end{equation} 
It follows from  \eqref{jac.ij} that for  $j \ne i$,  
 \begin{equation}\label{eq.hij} 
 h_{ij} :  =    \frac{\partial}{\partial x_j} S_i(x ^* , \boldsymbol{0})=
- \tfrac{1}{\beta} \frac{  \nabla^2_{ij} f (x ^* ) }{ 1+  \tfrac{\lambda}{\beta} \tfrac{r'(|x ^*_i |)}{|x ^*_i |}  }.
\end{equation} 
We then combine this with \eqref{jac.ieps}--\eqref{eq.hij} to derive  \eqref{jac.h}.  
 The boundedness of $\Lcal(F^0)$ trivially leads to  the Lipschitz continuity of $S$.

(ii) This is straightforward by applying  \cite[Lemma 2.13 and 2.14]{davis2022proximal}. The    Jacobian $DT(x ^* , \boldsymbol{0})$ is obtained directly by (i). 
\end{proof}

Now we show that every strict saddle point of \eqref{pp} is an unstable fixed point of $T$. 

\begin{theorem}(Unstable Fixed Point) \label{th.Ufp2}
Suppose Assumptions \ref{ass.1}-\ref{ass.4} hold.   
If $x^*$ is a strict saddle point of \eqref{pp}, then the corresponding $(x^*,\epsilon^*)$ with $\epsilon^*= \boldsymbol{0}$ is an unstable fixed point of $T$. 
Moreover, for $\alpha < \beta/\rho$,  $DT(x^*, \boldsymbol{0})$ is invertible. \label{ThLs2.l2}
\end{theorem} 
\begin{proof}
By Theorem \ref{ModelId2}(ii),  it suffices to check the eigenvalues of 
  \begin{equation} \label{eq.hiw}    (1-\alpha)I  + \alpha H_{\Ical^*\Ical^*}  =  I-  \tfrac{\alpha}{\beta} (I+ \tfrac{\lambda}{\beta}U^*_{\Ical^*\Ical^*})^{-1} \nabla^2_{\Ical^*\Ical^*} F(x^*),\end{equation}
  where we have $U^*_{ii} = u(x_i^*,0) < \infty, i\in \Ical^*$. 
If $\lambda_{\min}(\nabla^2_{\Ical^*\Ical^*} F(x^*)) < 0$, it follows from  Lemma \cite[Proposition 6.1]{serre2010matrices} that $ \lambda_{\min}( (I+ \tfrac{\lambda}{\beta}U^*_{\Ical^*\Ical^*})^{-1} \nabla^2_{\Ical^*\Ical^*} F(x^*))< 0$. 
 Hence, $DT(x^*, \boldsymbol{0})$ has an eigenvalue greater than 1, meaning $(x^*,  \boldsymbol{0})$ is the unstable fixed point of $T$.  
 Moreover,  since $\|   (I+ \tfrac{\lambda}{\beta}U^*_{\Ical^*\Ical^*})^{-1} \nabla^2_{\Ical^*\Ical^*} F(x^*)\| \le \|  (I+ \tfrac{\lambda}{\beta}U^*_{\Ical^*\Ical^*})^{-1}\| \| \nabla^2_{\Ical^*\Ical^*} F(x^*)\| \le \| \nabla^2_{\Ical^*\Ical^*} F(x^*)\|  \le \rho$  by \eqref{hess is bounded} and  
$ \nabla^2_{\mathcal{I}^*\Ical^*}  F(x^*)$ does not have 0 eigenvalue  by the strict saddle property,  we know by \eqref{hess is bounded} that 
$DT(x^*,\epsilon^*)$ is invertible   for  $\alpha < \beta/\rho$. 
\end{proof}

Now we can derive the global escape property for the ${\rm DIRL}_2$ algorithm.
 \begin{theorem}(Global Escape)
    Consider problem \eqref{pp} with strict saddle property and the ${\rm DIRL}_2$ algorithm with  Assumptions \ref{ass.1}-\ref{ass.4}. 
    Let  
     $\alpha < \beta/\rho$ and $\alpha \in(0, \frac{1}{1+L_S})$ where $L_S$ is as defined in Theorem \ref{ModelId2}. Let $\{(x^k,\epsilon^k)\}$ be the sequence generated by the Algorithm \ref{alg2}. Then, the set of initial points attracted by strict saddle points has 
    zero Lebesgue measure. 
\end{theorem}
\begin{proof}
    It follows directly from Theorems \ref{ModelId2} and \ref{th.Ufp2} that: (i) map $T$ is local $C^1$-smooth around any stationary point; (ii) all stationary points $x^*$ are unstable fixed points; (iii) $T$ is a lipeomorphism. Then by Corollary \ref{coro.saddle}, the set of initial conditions attracted by such fixed points has 
    zero Lebesgue measure. 
\end{proof}

\appendix
\section{Proof supplementary} 
\renewcommand{\thelemma}{A.\arabic{lemma}}
\renewcommand{\thetheorem}{A.\arabic{theorem}}
\renewcommand{\thedefinition}{A.\arabic{definition}}

\begin{lemma}\label{lem.nonexpansive}  Given $\hat w, \tilde w \in \mathbb{R}^n_+$ and $\hat z, \tilde z \in \mathbb{R}^n$, it holds that 
\[  \| \Sb_{\hat w}(\hat z) - \Sb_{\tilde w}(\tilde z)  \| \le \| \hat z - \tilde z\| + \| \hat w - \tilde w\|. \]
\end{lemma} 
\begin{proof} We consider the following three cases.

 (a) If $| \hat z_i | \le  \hat w_i$ and $|\tilde y_i| \le  \tilde w_i$, 
 then $\Sb_{  \hat w_i } ( \hat z_i ) = 0 $ and  $\Sb_{  \tilde w_i } ( \tilde  z_i ) = 0$. 
 In this case, $\Sb_{  \hat w_i } ( \hat z_i )  - \Sb_{  \tilde w_i } ( \tilde  z_i )  
 = |0-0| = 0\le | \hat z_i - \tilde z_i |$. 

 (b)  If $| \hat z_i | >  \hat w_i$ and $|\tilde z_i| \le  \tilde w_i$ (or vice versa),  
 then $\Sb_{ \hat w_i} ( \hat z_i ) = \sign(\hat z_i) (|\hat z_i|-  \hat w_i ) $ and 
  $\Sb_{ \tilde w_i} ( \tilde  z_i ) = 0$. 
  In this case,   $\Sb_{  \hat w_i} ( \hat z_i ) -\Sb_{ \tilde w_i} ( \tilde  z_i ) = 
| |\hat z_i|-  \hat w_i |$. 
It follows that   
\[\begin{aligned}
 | \hat z_i - \tilde z_i | \ge &\ | \hat z_i | - | \tilde z_i | \ge | \hat z_i | -   \tilde w_i \\
   =  & \ | | \hat z_i | -    \hat w_i  + \lambda ( \hat w_i -   \tilde w_i)/\beta|  \\
   \ge &\  | \Sb_{  \hat w_i } ( \hat z_i ) -\Sb_{  \tilde w_i } ( \tilde  z_i ) | - |\lambda ( \hat w_i -   \tilde w_i)/\beta|,
\end{aligned}\]
implying $| \Sb_{  \hat w_i } ( \hat z_i ) -\Sb_{  \tilde w_i } ( \tilde  z_i ) |  \le   | \hat z_i - \tilde z_i | 
+  | \hat w_i -   \tilde w_i| .$
 
 (c)  If $| \hat z_i | >  \hat w_i$ and $|\tilde z_i| >  \tilde w_i$, 
  then $| \Sb_{  \hat w_i } ( \hat z_i ) - \Sb_{  \tilde w_i } ( \tilde  z_i ) |
  =   | \sign(\hat z_i) (|\hat z_i|-  \hat w_i ) - \sign(\tilde z_i) (|\tilde z_i|-  \tilde w_i )|.$
If $\sign(\hat z_i) = \sign(\tilde  z_i)$,  
 $| \Sb_{  \hat w_i } ( \hat z_i ) - \Sb_{  \tilde w_i } ( \tilde  z_i ) |
  = |  (|\hat z_i|-  \hat w_i ) -  (|\tilde z_i|-  \tilde w_i )| \le 
   | \hat z_i - \tilde z_i | 
+  | \hat w_i -   \tilde w_i|$. If $\sign(\hat z_i) = -  \sign(\tilde  z_i)$,
$| \Sb_{  \hat w_i } ( \hat z_i ) - \Sb_{  \tilde w_i } ( \tilde  z_i ) | =
 (|\hat z_i|-  \hat w_i ) + (|\tilde z_i|-  \tilde w_i ) \le 
  |\hat z_i| + |\tilde z_i|  = |\hat z_i - \tilde z_i|.
$
   
In all the three cases,  
we show $| \Sb_{  \hat w_i } ( \hat z_i ) - \Sb_{  \tilde w_i } ( \tilde  z_i ) | \le |\hat z_i - \tilde z_i| 
+  | \hat w_i -   \tilde w_i|.$
Therefore,
$\| \Sb_{ \hat w}(\hat w)  - \Sb_{  \tilde w}(\tilde w) \| \le \| \hat z - \tilde z\| + 
 \| \hat w -   \tilde w\|.$
 \end{proof}

\newpage

\bibliographystyle{siamplain}
\bibliography{references}

\end{document}


\maketitle

\section{A detailed example}

Here we include some equations and theorem-like environments to show
how these are labeled in a supplement and can be referenced from the
main text.
Consider the following equation:
\begin{equation}
  \label{eq:suppa}
  a^2 + b^2 = c^2.
\end{equation}
You can also reference equations such as \cref{eq:matrices,eq:bb} 
from the main article in this supplement.

\lipsum[100-101]

\begin{theorem}
An example theorem.
\end{theorem}

\lipsum[102]
 
\begin{lemma}
An example lemma.
\end{lemma}

\lipsum[103-105]

Here is an example citation: \cite{KoMa14}.

\section[Proof of Thm]{Proof of \cref{thm:bigthm}}
\label{sec:proof}

\lipsum[106-112]

\section{Additional experimental results}
\Cref{tab:smfoo} shows additional
supporting evidence. 

\begin{table}[htbp]
\footnotesize
  \caption{Example table.}\label{tab:smfoo}
\begin{center}
  \begin{tabular}{|c|c|c|} \hline
   Species & \bf Mean & \bf Std.~Dev. \\ \hline
    1 & 3.4 & 1.2 \\
    2 & 5.4 & 0.6 \\ \hline
  \end{tabular}
\end{center}
\end{table}

\bibliographystyle{siamplain}
\bibliography{references}